\documentclass[a4paper,12pt]{article}
\usepackage{amssymb}
\usepackage{amsthm}
\usepackage{amsxtra}
\usepackage{amsmath, amscd}
\usepackage{amsfonts}
\usepackage{xr}
\usepackage{wasysym}
\usepackage[all,cmtip]{xy}

\begin{document}

\def\sect{\section}

\newtheorem{thm}{Theorem}[section]
\newtheorem{cor}[thm]{Corollary}
\newtheorem{lem}[thm]{Lemma}
\newtheorem{prop}[thm]{Proposition}
\newtheorem{propconstr}[thm]{Proposition-Construction}
\newtheorem{remcor}[thm]{Remark/Corollary}

\theoremstyle{definition}
\newtheorem{para}[thm]{}
\newtheorem{ax}[thm]{Axiom}
\newtheorem{conj}[thm]{Conjecture}
\newtheorem{defn}[thm]{Definition}
\newtheorem{notation}[thm]{Notation}
\newtheorem{rem}[thm]{Remarks}
\newtheorem{remark}[thm]{Remark}
\newtheorem{question}[thm]{Questions}
\newtheorem{example}[thm]{Example}
\newtheorem{problem}[thm]{Problem}
\newtheorem{excercise}[thm]{Exercise}
\newtheorem{ex}[thm]{Exercise}
\newtheorem{comments}[thm]{Comments}
\newtheorem{rems}[thm]{Remarks}

\def\Bbb{\mathbb}
\def\cal{\mathcal}
\def\mL{{\mathcal L}}
\def\mC{{\mathcal C}}

\overfullrule=0pt

\def\si{\sigma}
\def\prf{\smallskip\noindent{\it        Proof}. }
\def\call{{\cal L}}
\def\nat{{\Bbb  N}}
\def\la{\langle}
\def\ra{\rangle}
\def\inv{^{-1}}
\def\ld{{\rm    ld}}
\def\trdeg{{tr.deg}}
\def\dim{{\rm   dim}}
\def\th{{\rm    Th}}
\def\rest{{\lower       .25     em      \hbox{$\vert$}}}
\def\restr#1{{\lower       .25     em      \hbox{$\vert\scriptstyle{#1}$}}}
\def\ch{{\rm    char}}
\def\zee{{\Bbb  Z}}
\def\conc{^\frown}
\def\acl{{\rm acl}}
\def\cls{cl}
\def\cals{{\cal S}}
\def\mult{{\rm  Mult}}
\def\calv{{\cal V}}
\def\aut{{\rm   Aut}}
\def\ffi{{\Bbb  F}}
\def\ffiti{\tilde{\Bbb          F}}
\def\degs{deg_\si}
\def\calx{{\cal X}}
\def\gal{{\cal G}al}
\def\cl{{\rm cl}}
\def\loc{{\rm locus}}
\def\calg{{\cal G}}
\def\calq{{\cal Q}}
\def\calr{{\cal R}}
\def\caly{{\cal Y}}
\def\aff{{\Bbb A}}
\def\cali{{\cal I}}
\def\calu{{\cal U}}
\def\epsilon{\varepsilon} 
\def\Uu{{\cal U}}
\def\rat{{\Bbb Q}}
\def\ga{{\Bbb G}_a}
\def\gm{{\Bbb G}_m}
\def\cee{{\Bbb C}}
\def\ree{{\Bbb R}}
\def\frob{{\rm Frob}}
\def\Frob{{\rm Frob}}
\def\fix{{\rm Fix}}
\def\Uu{{\cal U}}
\def\proj{{\Bbb P}}
\def\sym{{\rm Sym}}
 
\def\dcl{{\rm dcl}}
\def\calm{{\mathcal M}}
\def\calp{{\cal P}}

\font\helpp=cmsy5
\def\semdp
{\hbox{$\times\kern-.23em\lower-.1em\hbox{\helpp\char'152}$}\,}

\def\dnfo{\,\raise.2em\hbox{$\,\mathrel|\kern-.9em\lower.35em\hbox{$\smile$}
$}}
\def\dnf#1{\lower1em\hbox{$\buildrel\dnfo\over{\scriptstyle #1}$}}
\def\dfo{\;\raise.2em\hbox{$\mathrel|\kern-.9em\lower.35em\hbox{$\smile$}
\kern-.7em\hbox{\char'57}$}\;}
\def\df#1{\lower1em\hbox{$\buildrel\dfo\over{\scriptstyle #1}$}}        
\def\stab{{\rm Stab}}
\def\qfcb{\hbox{qf-Cb}}
\def\perf{^{\rm perf}}
\def\sipm{\si^{\pm 1}}

\def\calc{{\cal C}}
\def\ione{{\it (I1)}}
\def\itwo{{\it (I2)}}
\def\itwow{{\it (I2w)}}
\def\ithreew{{\it (I3w)}}
\def\ithree{{\it (I3)}}
\def\ifour{{\it (I4)}}
\def\ifourw{{\it (I4w)}}
\def\bigudot{\hbox{$\bigcup\kern-.75em\cdot\;\,$}}
\def\caln{{\cal N}}
\def\cald{{\cal D}}
\def\calm{{\cal M}}
\def\calu{{\cal U}}
\def\Div{\hbox{Div}}
\def\pic{\hbox{Pic}}
\def\G{{\cal G}}
\def\sime{{\hbox{$\kern-.19em\sim$}}}
\def\Th{{\rm Th}}
\def\dnfo{\,\raise.2em\hbox{$\,\mathrel|\kern-.9em\lower.35em\hbox{$\smile$}
$}}
\def\dfo{\;\raise.2em\hbox{$\mathrel|\kern-.9em\lower.35em\hbox{$\smile$}
\kern-.7em\hbox{\char'57}$}\;}
\def\vlabel{\label}
\def\SCF{{\rm SCF}}

\title{Amalgamation of types in
pseudo-algebraically closed fields and applications}

\author{Zo\'e Chatzidakis\thanks{Most of this work was done while the
    author was partially supported by
    MRTN-CT-2004-512234 and by ANR-06-BLAN-0183, while a member of the
    Equipe de Logique Math\'ematiques (UMR 7056), in University Paris~Diderot. The  work  achieved its
    final form in 2017, while  the author was partially supported by
    ANR-13-BS01-0006.}\\ DMA (UMR 8553), Ecole Normale Sup\'erieure\\
    CNRS, PSL Research University
}
\maketitle

\begin{abstract} This paper studies unbounded PAC fields and shows an
  amalgamation result for types over algebraically closed sets. It
  discusses various applications, for instance that omega-free PAC
  fields have the property NSOP3. It also contains a description of
  imaginaries in PAC fields.  
\end{abstract}

\section*
{\bf Introduction} Pseudo-algebraically closed fields (henceforth abbreviated
by PAC) were introduced by
Ax in his famous paper \cite{Ax} on the theory of finite fields. The
elementary theory of arbitrary PAC fields, studied among others by Cherlin-Van den
Dries-Macintyre \cite{[CDM]} and by Ershov \cite{Er}, puts in light an interesting dichotomy:
definable sets are given, on the one hand by classical algebraic data, and on
the other hand by elementary statements concerning the Galois
group. Many of the properties of the theory of a PAC field thus reduce
to the corresponding properties of its Galois group. For instance, if
the subfield of algebraic numbers of the PAC field $F$ is decidable,
then ${\rm Th}(F)$ will be decidable if and only if the ``theory'' of
its absolute Galois group is decidable. One also knows that the
structure of their models 
is complicated: a result of Duret (\cite{Du}) asserts that a PAC field
which is not separably closed has the independence property.

Interest for the model theory of PAC fields revived in the mid 90's, when
 Hrushovski and Pillay (\cite{HP}) were able to use stability
theoretic techniques for groups definable in pseudo-finite fields, and
more generally in bounded PAC fields (a field  is {\em bounded} if for each
$n>1$ it has only finitely many algebraic  extensions of degree $n$). It was then
observed that bounded PAC fields have a simple theory, because
they satisfy the independence theorem (1991 result of Hrushovski, only
published in 2005, \cite{H}).  Other results with a stability-theoretic flavour
followed: in \cite{[C1]}, the author shows that a PAC field with a
simple theory is necessarily bounded; a weak notion of independence is
defined, and shown to be implied (in any field) by non-forking. In \cite{[C2]}, the
study of unbounded PAC fields is continued, with emphasis on the theory
of $\omega$-free PAC fields. The author shows that for these fields,
forking is the transitive closure of  weak independence, and shows  versions of
the independence theorem for various independence notions, the most difficult one being that $\omega$-free
PAC fields of chararacteristic $0$ satisfy the independence theorem with
independence being the genuine non-forking. This last
result is quite 
surprising, given that the theories of $\omega$-free PAC fields are not
simple. This suggested that more can be done on unbounded PAC fields,
and that  their study might
provide an insight of good behaviours of models of non-simple theories.

In this paper we continue the investigation of the behaviour of unbounded PAC fields.  Our main result is an amalgamation
result for types, similar to the (weak) independence theorem of
\cite{[C2]}.  This result 
(Theorem \ref{thm1}) isolates the conditions under which amalgamation of types is
possible.  It has various consequences, notably a weak independence
theorem over models for PAC fields $F$ such that $S\G(F)$ has a simple
theory 
(Theorem \ref{thm25}),
and the fact that Frobenius fields satisfy NSOP$_3$ (see \ref{thm41} and
\ref{cor42}).  It  also
appears as an ingredient in the description of imaginaries in PAC
fields of finite degree of imperfection: an imaginary of the
PAC field $F$ is equi-definable with a finite collection of pairs
$(a,D)$, where $a$ is a tuple of elements of $F$ and $D$ is an
imaginary of $S\G(F)$ (Theorem \ref{thm31}). We show by an example that
this result is best possible. 

The hope that PAC fields might provide good examples of things happening
beyond simplicity was vindicated. Recent results of Chernikov and
Ramsey (\cite{CR}, Theorem~6.2) show that the weak independence theorem proved in \cite{[C2]} for
Frobenius fields implies that the theory of a Frobenius field is
NSOP$_1$. Thus these fields provide a large family of new examples of structures
with an NSOP$_1$ theory. This is particularly useful as very few examples of theories with
NSOP$_1$ were known.  It can be hoped that a further study of
these PAC fields might lead to new insight on NSOP$_1$ theories. The $\omega$-free PAC fields are particularly nice Frobenius
fields, in which types and definable sets are well understood. As we  show here, imaginaries are
equally well understood. 

Clearly, the
connections between the neo-stability properties of the Galois group of
a PAC field 
and those of the field also need to be explored further. Results of Nick
Ramsey (\cite{R}) suggest this is the case for the properties NSOP$_1$
and NTP$_1$. \\

The paper is organised as follows. In section 1, after setting
up the notation, we recall or prove 
 some  technical results on fields and profinite groups. Section 2  contains the main result of
this paper, Theorem \ref{thm1}, as well as various independence theorems
and SOP$_n$ properties for $n\geq 3$. We  conclude  section 2 with
some questions. Section 3 develops the part of the
logic of complete systems which is interpretable in fields. In
particular, it sets up the  formalism which will enable us to deal with definable
sets. This is applied in section 4 (Theorem \ref{thm31}) to give the description of imaginaries
of PAC fields $F$ of finite degree of imperfection.

\section{Notation and preliminary results}

Recall first 
  that a field $F$ is PAC if every absolutely irreducible variety
  defined over $F$ has an $F$-rational point. Equivalently, if $F$ is
  existentially closed in any regular extension. 
In this section we set up the notation,  recall some classical 
results on PAC fields, and give two additional lemmas. We assume
familiarity with 
elementary results on field extensions, see e.g. Chapter III of \cite{[L]}.

\para\vlabel{not}{\bf Notation, conventions}. We work in the usual
language of rings ($\{+,-,\cdot,0,1\})$, sometimes expanded by adding
constants for a $p$-basis. The separable closure of a
field $K$ is denoted by $K^s$, and its absolute Galois group $\gal(K^s/K)$
by $\G(K)$. If $A\subseteq K$, then $\acl(A)$ denotes the model-theoretic
closure of $A$ in the sense of $\Th(K)$. It is known that $K$ is a regular
extension of $\acl(A)$. 

We will often work inside the separable closure of a field $K$.  In
that case, we will denote by SCF the theory $\Th(K^s)$, the notation
$tp_{\SCF}(\ \ )$ will refer to the type in the field $K^s$.  We use
the notation $\acl_{K^s}(A)$ to denote the algebraic closure in the
sense of $\Th(K^s)$, i.e., the smallest subfield of $K^s$ containing
$A$ and of which $K^s$ is a regular extension.  We will say that two
subsets of $K$ (or of $K^s$) are SCF-independent over some $E$ if they
are independent in the sense of $\Th(K^s)$.  

In addition, unless otherwise specified, all fields will be subfields of
some large algebraically closed field $\Omega$. If $A,B$ are two
subfields, then $AB$ denotes the composite field.

\bigskip
An extremely useful and fundamental result on PAC fields is the so-called
``embedding lemma'' of Jarden and Kiehne:

\begin{thm}\vlabel{pac4} {\rm (Lemma 20.2.2 in \cite{[FJ]})} 
Let $E/L$ and $F/M$ be separable field extensions satisfying:
$E$ is countable  and  
$F$ is an 
  $\aleph_1$-saturated PAC field; if ${\rm char}(F)=p>0$, assume in
  addition that $[E:E^p]\leq [F:F^p]$. Assume that there is an isomorphism
  $\varphi_0:L^s\to M^s$ such that $\varphi_0(L)=M$, and a commutative
  diagramme
\begin{equation*}
\begin{CD} \G(E) @<\Phi<< \G(F) \\ 
@V{\rm res}VV @VV{\rm res}V \\ 
\G(L) @<\Phi_0<< \G(M)
\end{CD}
\end{equation*}
where $\Phi_0:\si\mapsto \varphi_0\inv\si\varphi_0$, is the {\em dual} of $\varphi_0$, and $\Phi$ is a (continuous)
homomorphism. Then $\varphi_0$ extends to an embedding $\varphi:E^s\to
F^s$, with dual $\Phi$, and such that $F/\varphi(E)$ is separable. 
\end{thm}

\begin{rem} \vlabel{rempac4} We will  use the following essentially immediate
  consequences of this result. 
\begin{enumerate}
\item We may replace the countability hypothesis on $E$ by asking $F$ to
  be $|E|^+$-saturated. The proof is identical.
\item We will usually have that the extensions $E/L$ and $F/M$ are
  regular. This means that the restriction maps $\G(E)\to \G(L)$ and
  $\G(F)\to \G(M)$ are onto. Note that the conclusion will then be that
  $F/\varphi(L)$ is regular. Similarly, if $\Phi$ is onto, then the
  extension $F/\varphi(E)$ will be regular. 
\item (Notation as above.) Let $E'$ be a Galois extension of $E$
  containing $L^s$, and  $\Phi':\G(F)\to
  \gal(E'/E)$ such that the following diagramme commutes: 
\begin{equation*}
\begin{CD} \gal(E'/E) @<\Phi'<< \G(F) \\ 
@V{\rm res}VV @VV{\rm res}V \\ 
\G(L) @<\Phi_0<< \G(M)
\end{CD}
\end{equation*}
As $\G(F)$ is projective, the map $\Phi'$ factors through a homomorphism
$\Phi:\G(F)\to \G(E)$ (see Theorem 11.6.2  in \cite{[FJ]}). Applying the
embedding lemma therefore gives us an embedding $\varphi':E'\to F^s$,
with dual $\Phi'$.  

\end{enumerate}
\end{rem}

\subsection*{Complete systems associated to profinite groups}
Cherlin, Van den
Dries and Macintyre show in \cite{[CDM]} how to associate to any profinite group $G$
a  structure $SG$ in an $\omega$-sorted language $\call_G$, called the
{\it complete system
of $G$}, which encodes precisely the
inverse system of all finite continuous quotients of $G$. The
functor $G\mapsto SG$ is a contravariant functor, and defines a duality
between  the category
of profinite groups with continuous epimorphisms and the category of
complete systems with embeddings. The functor dual to $S$ is the functor
$G$ which to a complete
system $S$ associates the inverse limit of the inverse system of finite
groups given by $S$. An important remark, which is at the core of the results of Cherlin Van
den Dries and Macintyre, is that the functor $S\G$ commutes with
ultraproducts and therefore with ultrapowers: If $\calu$ is an ultrafilter on a set $I$ and $K$ is a field, then
$S\G(K^\calu)\simeq (S\G(K))^\calu$, where the second ultraproduct is
taken in the $\omega$-sorted context (i.e., sort by sort). Hence, $K\equiv L$ implies $S\G(K)\equiv
S\G(L)$. In an unpublished manuscript, Cherlin, Van den Dries and
Macintyre  also show that this $\omega$-sorted  logic on $S\G(K)$ is
in some sense the strongest logic of the Galois group $\G(K)$ which is
interpretable in the field $K$. For more details on
complete systems and their logic,
see \cite{[CDM]} or the Appendix of \cite{[C2]}. \\
We will first briefly recall
the notation and definitions for arbitrary profinite groups, before
going to the setting of Galois groups.

\para\vlabel{cs1} {\bf Definition of the complete system of a profinite group}. \\
Let $G$ be a profinite group, and $\call_G$  be the
$\omega$-sorted language with sorts indexed by the positive integers, and
with non-logical symbols $\{\leq,
C,P,1\}$, where $\leq$ and $C$ are binary relations,  $P$ is a
ternary relation and $1$ is a constant symbol.  The complete system
associated to $G$ is the $\call_G$-structure $S(G)$, with universe the
disjoint union $\bigudot_NG/N$ where $N$ ranges over all normal open
subgroups of $G$. An element of $G/N$, i.e. a coset $gN$, will be of
sort $n$ if and only if $[G:N]\leq n$, and $1=G$ is the only element
of sort $1$. We have $gN\leq hM\iff N\subseteq M$, $C(gN,hM)\iff
gN\subseteq hM$,
and $P(g_1N_1,g_2N_2,g_3N_3)\iff N_1=N_2=N_3$ and $g_1g_2N_1=g_3N_1$.
The class of complete systems of profinite groups is the class of models
of  a theory $T_G$. The functor $S$ defines a duality between the
category of profinite groups with continuous epimorphisms and the
category of models of $T_G$ with embeddings.

\para\vlabel{cs2}{\bf Complete systems of Galois groups, subsystems,
double duals}. \\ Let $F$ be a
field, $E$ a Galois extension of $F$, and $G=\gal(E/F)$. The universe of
$SG$ is the disjoint union of all $\gal(L/F)$ where $L$ is a finite
Galois extension of $F$ contained in $E$. The elements of sort $n$ with
be the Galois 
groups of size $\leq n$. The language $\call_G$ is interpreted as
follows:  $\gal(F/F)=1$; whenever
$L_1\supseteq L_2$, $C\cap (\gal(L_1/F)\times \gal(L_2/F))$ is the graph
of the restriction maps $\gal(L_1/F)\to \gal(L_2/F)$  and $\leq$ contains
$\gal(L_1/F)\times \gal(L_2/F)$; 
 the ternary relation $P$ encodes the graph of multiplication on each
$\gal(L/F)$. 

\smallskip
A subset  $S$ of $SG$ is a {\it subsystem}\footnote{Warning: in earlier papers by the author
 they are called {\em substructures}.} of $SG$ if it has the
following two properties: (i) 
$\forall \si,\tau\in S,\, \exists \rho\in S\ (\rho \leq \si \land \rho \leq
\tau)$; (ii) If $\si\in S$ and $\tau\geq \si$, then $\tau\in S$. If
$A\subset SG$, then $\langle A\rangle$ denotes the subsystem of $SG$ generated
by $A$. 

One sees easily that if $S$ is a subsystem of $SG$, then
$S=S\gal(M/F)$, where $M$ is the composite of all Galois extensions $L$
such that $S$ contains $\gal(L/F)$. The inclusion map $S\subset SG$ and
the restriction map $\gal(E/F)\to \gal(M/F)$ are {\it dual of each other}.

Let $F_1$ and $F_2$ be fields, and $\varphi:F_1^s\to F_2^s$ an
embedding such that $\varphi(F_1^s)\cap F_2=F_1$. We then get a
continuous epimorphism $\Phi:\G(F_2)\to \G(F_1)$, defined by $\si\mapsto
\varphi\inv\si\varphi$ (the {\em dual} of $\varphi$). Applying the functor $S$  to $\Phi$ gives us
an embedding $S\G(F_1)\to S\G(F_2)$, defined as follows: if $L_1$ is
a finite Galois extension of $F_1$ and $\si\in \gal(L_1/F_1)$, then
$S\Phi(\si)$ is the unique element of $\gal(F_2\varphi(L_1)/F_2)$ extending
the element $\varphi\si\varphi\inv$ of $\gal(\varphi(L_1)/\varphi(F_1))$.
We  call the map $S\Phi$ the {\it double dual of $\varphi$}.

\begin{thm}\vlabel{cdm2}
 (Cherlin, Van den Dries, Macintyre \cite{[CDM]}).
Let $F_1$ and $F_2$ be PAC fields, separable over
a common subfield $E$. The following conditions are equivalent:
\begin{enumerate}
\item[(1)]{$F_1\equiv_E F_2$.}
\item[(2)]
\begin{enumerate}\item[(i)] $F_1$ and $F_2$ have the same degree of
imperfection,
\item[(ii)]{There is  $\varphi\in \G(E)$ such that
$\varphi(F_1\cap E^s)=F_2\cap E^s$, and the double dual
$S\Phi:S\G(F_1\cap
E^s)\to S\G(F_2\cap E^s)$  of $\varphi$, is a partial elementary
$\call_G$-map $S\G(F_1)\to S\G(F_2)$.  (In particular, $S\G(F_1)\equiv S\G(F_2)$).}
\end{enumerate}
\end{enumerate}
\end{thm}

 From this result, one easily deduces a description of types:
\begin{thm}\vlabel{cdm3} Let $F$ be a PAC field,
separable over some subfield $E$.
Let $a$ and $b$ be tuples of elements of $F$, and $A=\acl_{K^s}(E,a)\cap
F$,
$B=\acl_{K^s}(E,b)\cap F$. The following conditions are equivalent:
\begin{enumerate}
\item{$tp(a/E)=tp(b/E)$.}
\item{There is an $E$-isomorphism $\varphi:A^s\to B^s$, with
$\varphi(a)=b$, $\varphi(A)=B$, such that the double dual $S\Phi:S\G(A)
\to S\G(B)$ is a
partial elementary $\call_G$-map of $S\G(F)$.}
\end{enumerate}
\end{thm}

\para \vlabel{acl}{\bf Important facts and remarks}. If $F$ is a PAC field and
$A\subset F$, then $\acl(A)=\acl_{F^s}(A)\cap F$ (see 4.5 in \cite{[CP]}). Let
$E\subset A,B$ be subfields of $F$, and assume  that $A$ and $B$ are
SCF-independent over $E$. If $char(F)=p>0$ and $[F:F^p]<\infty$ assume
moreover that $E$
contains a $p$-basis of $F$. Then
$\acl_{F^s}(AB)=(\acl_{F^s}(A)\acl_{F^s}(B))^s$. Hence we also have
$\acl(AB)=(\acl(A)\acl(B))^s\cap F$. 

\subsection*{Frobenius  and $\omega$-free PAC fields}
\begin{defn}\vlabel{free1}\vlabel{frob1} \begin{enumerate}
\item A profinite group $G$ has the {\em  embedding property}  if for
  any finite groups $A$, $B$, 
  whenever $f:G\to A$ and $g:B\to A$, $f':G\to B$ are (continuous)
  epimorphisms, then there exists an epimorphism $h:G\to B$ such that
  $f=g\circ h$:

$$
\xymatrix{ & G \ar[d]^{f} \ar@{-->}[dl]_-{\exists h} \\
 B \ar[r]^{g}& A
}
$$
This property translates into a property of $SG$ which is axiomatisable in the language
  $\call_G$. See section
  24.3 of \cite{[FJ]} for more details and properties of these groups. 
\item A {\em Frobenius field} is a PAC
field whose absolute Galois group $\G(F)$ has the  embedding property. 

\item Recall that a PAC
field $F$ is {\em $\omega$-free} if whenever $F_0\prec F$ is countable, then
$\G(F_0)\simeq \hat F_\omega$, the free profinite group on $\aleph_0$
generators. In particular (or equivalently, using a result of Iwasawa), all finite groups occur as finite quotients
of $\G(F)$, and $F$ is Frobenius. 
\end{enumerate}\end{defn}
\noindent
Being Frobenius is an elementary
property of a field $F$. When dualized,  and if $S\G(F)$ is countable, 
it says that any $\call_G$-isomorphism between two finite subsystems
of $S\G(F)$ extends to an automorphism of $S\G(F)$.  In particular this
implies the following:

\begin{thm}
(\cite{[CDM]}) If $F$ is a Frobenius field, then any
$\call_G$-isomorphism between two 
subsystems of $S\G(F)$ is elementary. 
\end{thm}

\noindent
Thus, in  Theorems \ref{cdm2} and \ref{cdm3}, the conditions stating
that the partial maps $S\Phi$ are elementary can be removed.

\para $\omega$-sorted logic behaves very much like ordinary one-sorted
logic, provided one works sort by sort. Our $\omega$-sorted structure
$S$ can be viewed as the countable union of structures $S_n$, $n\geq 1$,
where 
each $S_n$ has 
universe the elements of sort $\leq n$, and is a structure in
the language with $n$ sorts, relational symbols $P,C$ and $\leq$,
constant symbol $1$. The theory of $S$ is then naturally the limit of
the theories of the $S_n$'s. For instance, let $G$ be a profinite group
with the embedding property, $SG$ its complete system. Then the above
characterisation of countable models of ${\rm Th}(SG)$ translates into: ${\rm Th}(SG)$ is $\aleph_0$-categorical (see
\cite{[CDM]}). 
Notions such as stability or
$\omega$-stability easily generalise: one just counts types in each
sort. Notions which are local immediately generalise, as they only
involve finitely many sorts. For instance, the usual definition of
forking of  a
formula over a set; and therefore forking of a type: a type will fork
over a set if it
contains a formula which forks over that set.  Hence one can define the
property of a theory of being simple. Note that $\Th(S)$ is simple if
and only if $\Th(S_n)$ is simple for every $n\geq 1$. We will use the fact that
the results of 
Kim and Pillay characterizing simple theories via the properties of the
forking relation go through.

\bigskip\noindent
Before proving our two technical lemmas, we first recall some results
from \cite{[C2]}:

\begin{lem}\vlabel{lem1} {\rm (Lemma 2.1 in \cite{[C2]})} Let $A$, $B$ be fields (contained in 
$\Omega$), and assume that the field composite $AB$ is a regular extension of $A$ and of $B$.
If $E=A\cap B$, then $E^s=A^s\cap B^s$.
\end{lem}

\begin{lem}\vlabel{lem5} Let $A$, $B$, $C$, $E$
be separably 
closed fields ($\subset \Omega$), with
$A,B,C$ separable extensions of $E$, and $AB$ a separable extension
of $A$ and of $B$, free from $C$ over $E$, and with $A\cap B=E$. Then

\begin{itemize}
\item[(i)] (a) $(AB)^s\cap (AC)^s(BC)^s=AB$; (b) $(AB)^sC\cap (AC)^s(BC)^s=ABC$
\item[(ii)] (a) $(AC)^s\cap (AB)^s(BC)^s=AC$; (b) $(AC)^sB\cap
  (AB)^s(BC)^s=ABC$.
\end{itemize}
\end{lem}

\prf Items (4) and (2) of Lemma 2.5  in
\cite{[C2]} give (i)(a)(b) and (ii)(a). By 
 (3) of that same lemma, we have $(AC)^s(AB)^s\cap
(BC)^s(AB)^s=C(AB)^s$, which implies $(AC)^sB\cap (AB)^s(BC)^s\subseteq
C(AB)^s\cap (AC)^sB\subseteq ABC$ by (i)(b), and gives us (ii)(b).

\begin{lem}\vlabel{lem6} Let $A$, $B$, $C$ (contained in $\Omega$) be regular extensions
of a field $E$, and assume that $AB$ is a regular extension of $A$ and
of $B$, that $A\cap B=E$, and that $AB$ is free from $C$ over $E$.
Consider the map
$$\rho: \gal((AB)^s(AC)^s(BC)^s/ABC)\to \G(AB)\times \G(AC)\times
\G(BC)$$ defined by 
$$\si\mapsto (\si\restr{(AB)^s},\si\restr{(AC)^s},\si\restr{(BC)^s}).$$ Then
the image of
$\rho$ is the subgroup of $\G(AB)\times \G(AC)\times \G(BC)$ consisting of
the triples $(\si_1,\si_2,\si_3)$ such that
$$\si_1\restr{A^s}=\si_2\restr{A^s},\ \si_1\restr{B^s}=\si_3\restr{B^s},\
\si_2\restr{C^s}=\si_3\restr{C^s}.$$
\end{lem}

\prf The compatibility conditions are clearly necessary, it remains to
show that they 
are sufficient. Let $(\si_1,\si_2,\si_3)\in\G(AB)\times \G(AC)\times
\G(BC)$ satisfy the required conditions. We will first show that there is some
$\si\in\gal(A^sB^sC^s/ABC)$ which agrees with $\si_1$ on $A^sB^s$, with $\si_2$
on $A^sC^s$ and with $\si_3$ on $B^sC^s$. First note that $\si_1$,
$\si_2$ and $\si_3$ all agree on $E^s$.

As $C$ is free from $AB$ over $E$, and is a regular extension of $E$, we
know that $C^s$ is linearly disjoint from $(AB)^s$ over $E^s$, and
therefore from $A^sB^sC$ over $CE^s$. Hence there is
$\si\in\gal(A^sB^sC^s/ABC)$ which agrees with $\si_1$ on $A^sB^s$ and
with $\si_2$ on $C^s$; by hypothesis, $\si$ therefore agrees with $\si_2$
on $A^sC^s$, and with $\si_3$ on $B^s$ and on $C^s$, i.e., on $B^sC^s$. 

By Lemma \ref{lem1}, $A^s\cap B^s=E^s$ and we may apply Lemma
\ref{lem5} to obtain:
\begin{itemize}
\item[$\bullet$]
the Galois 
extensions $(AB)^sC^s$, $(AC)^sB^s$ and $A^s(BC)^s$ are linearly disjoint over
$A^sB^sC^s$ (use (i)(b) and (ii)(b)), whence
\begin{multline*}\gal((AB)^s(AC)^s(BC)^s/A^sB^sC^s)\simeq \\
\gal((AB)^sC^s/A^sB^sC^s)\times \gal((AC)^sB^s/A^sB^sC^s)\times \gal(A^s(BC)^s/A^sB^sC^s);\end{multline*}
\item[$\bullet$] $A^sB^sC^s$ is a regular extension of $A^sB^s$, because
  $A^sB^sC^s\cap (AB)^s=A^sB^s$ (by (i)(a)), and therefore
 $$\gal((AB)^sC^s/A^sB^sC^s)\simeq \G(A^sB^s);$$
\item[$\bullet$] $A^sB^sC^s$ is a regular extension of $A^sC^s$, because
  $A^sB^sC^s\cap (AC)^s=A^sC^s$ (by (ii)(a)), so that $$\gal((AC)^sB^s/A^sB^sC^s)\simeq \G(A^sC^s);$$
\item[$\bullet$] $A^sB^sC^s$ is a regular extension of $B^sC^s$ (as
  above using (ii)(a)), 
so that $$\gal((BC)^sA^s/A^sB^sC^s)\simeq \G(B^sC^s).$$ 
\end{itemize}

\smallskip\noindent
This gives in particular that 
$$\gal((AB)^s(AC)^s(BC)^s/A^sB^sC^s)\simeq\G(A^sB^s)\times
\G(A^sC^s)\times \G(B^sC^s).$$
Hence, the automorphism $\si\in\gal(A^sB^sC^s/ABC)$ can be lifted
uniquely to an
element $\si$ of $\gal((AB)^s(AC)^s(BC)^s/ABC)$ which agrees with
$\si_1$ on $(AB)^s$, with $\si_2$ on $(AC)^s$ and with $\si_3$ on
$(BC)^s$.

\begin{lem}\vlabel{lem32} Let $E\subset A$ and $E_1$ be
algebraically closed subsets of a PAC field $F$, and assume that
$\varphi_0:E^s\to E_1^s$ is an isomorphism and  restricts to an
elementary map $E\to E_1$. If the characteristic is $p>0$ and
$[F:F^p]<\infty$, then assume that $E$ contains a $p$-basis of $F$. Let 
$S\Psi:S\G(A)\to S\subset S\G(F)$ be a partial  $\call_G$-elementary isomorphism extending the
double dual $S\Phi_0:S\G(E)\to S\G(E_1)$ of $\varphi_0$. Then in some
elementary extension $F^*$ of $F$, there is $B$, which is
SCF-independent from $F$ over $E_1$, and an
isomorphism $\varphi:A^s\to B^s$ sending $A$ to $B$,  extending $\varphi_0$ and with double dual
$S\Psi$. Moreover, $(B,E_1)$ realises $tp(A,E)$. 
\end{lem}

\prf Let $M$ be the Galois extension of $F$ with Galois group over $F$
corresponding to $S$, i.e., the restriction map ${\rm res}:\G(F)\to \gal(M/F)$ is
dual to $S\subset S\G(F)$. (Without the requirement that $B$ be
SCF-independent from $F$ over $E_1$, we could just apply Theorem
\ref{pac4}.) 

Choose any extension $\varphi$ of $\varphi_0$
to $A^s$ such that $\varphi(A^s)$ and $F^s$ are linearly disjoint over
$E_1^s$, and let $B=\varphi(A)$. 
Then the double dual $S\Phi$ of $\varphi$
extends $S\Phi_0$, and the dual $\Phi$ of $\varphi$ defines an
isomorphism  $\G(B)\to \G(A)$ which induces the dual $\Phi_0$ of $\varphi_0$,
$\Phi_0:\G(E_1)\to \G(E)$. 
 The dual $\Psi$ of $S\Psi$ defines an
isomorphism $\Psi:\gal(M/F)\to \G(A)$, which also induces $\Phi_0$ on $\G(E)$.
Consider the profinite group $$H=\{(\Phi\inv(\si),\Psi\inv(\si))\mid
\si\in\G(A)\}\subseteq \G(B)\times \gal(M/F).$$ Then $H$ is the graph of $\Psi\inv \Phi:\G(B)\to
\gal(M/F)$, and can be identified with a closed subgroup of
$\gal(B^sM/BF)\simeq \G(B)\times _{\G(E_1)}\gal(M/F)$. Let $L$ be the
subfield of $B^sM$ fixed by the elements of $H$. Since $H$ projects
onto $\G(B)$ and onto $G(S)=\gal(M/F)$, it follows that $L$ is a
regular extension of $B$ and of $F$, and that $\gal(B^sM/L)=H$
canonically identifies with $\G(B)$ and with $\gal(M/F)=G(S)$ via the
restriction maps. It follows that the restriction 
$\gal(B^sF^s/L)\to \G(F)$ is an isomorphism. By Theorem
\ref{pac4}, there is an elementary extension $F^*$ of $F$ containing $B$,
and 
such that $F^*\cap B^sF^s=L$. Then  the map $\varphi:A^s\to B^s$ is our
desired map: by construction, inside $S\G(F^*)$, we have $S\G(B)=S$,
and the double dual of $\varphi$ coincides with $S\Psi$, which is
an elementary map. This proves the first assertion, and \ref{cdm3} gives the moreover part. 

\begin{remark} \vlabel{remlem32} Let $E,E_1,A,\varphi_0$ be as above. Let $L$ be a Galois extension of $A$,  and assume that we have
  a partial elementary $\call_G$-map $S\Psi':S\gal(L/A)\to S'\subset
  S\G(F)$ which extends the restriction of $S\Phi_0$ to $S\gal(L\cap
  E^s/E)$. Then there is an elementary extension $F^*$ of $F$ such that the map
  $S\Psi'$ extends to an  $\call_G$-embedding $S\G(A)\to S\G(F^*)$. Thus in
  the above lemma, we may replace $S\G(A)$ by a subsystem. 
\end{remark}   

\begin{remark} Notation as in \ref{lem32}. The additional hypothesis in Lemma \ref{lem32} when $p>0$
and $[E:E^p]<\infty$ is actually not necessary but the proof needs to be
slightly modified. Indeed, the problem occurs if $[E:E^p]<[A:A^p]$ since
then one will get $[L:L^p]> [F:F^p]$ and we cannot apply the embedding lemma. This is not hard to fix: one selects
a $p$-basis $S$ of $B$ over $E$ (and therefore of $L$ over $F$), and forces
it to realise the generic $|S|$-type of $\Th(F^s)$ over $F^s$. Details
are left to the reader. 

\end{remark}

 \section{Amalgamation of types}
\begin{thm}\vlabel{thm1} Let $F$ be a PAC field, and let
$E,A,B,C_1,C_2$ be
algebraically closed subsets of $F$, with $E$ contained in $A,B,C_1,C_2$.
Assume that $A\cap B=E$, that  $A$ and $C_1$, and $B$ and $C_2$, are
SCF-independent over $E$, and that if the degree of imperfection of $F$ is
finite, then $E$ contains a $p$-basis of $F$. Moreover, assume that there is an
$E^s$-isomorphism $\varphi:C_1^s\to C_2^s$ such that $\varphi(C_1)=C_2$, and
that there is $S_0\subset S\G(F)$, and elementary (in $S\G(F)$)
isomorphisms
$$\displaylines{S\Psi_1:\langle S\G(C_1),S\G(A)\rangle \to \langle
S_0,S\G(A)\rangle\cr
S\Psi_2:\langle S\G(C_2),S\G(B)\rangle \to \langle
S_0,S\G(B)\rangle}$$
such that
\begin{enumerate}
\item[(i)]{$S\Psi_1$ is the identity on $S\G(A)$, $S\Psi_2$ is the identity
on $S\G(B)$, $S\Psi_i(S\G(C_i))=S_0$ and}
\item[(ii)]{if $S\Phi:S\G(C_1)\to S\G(C_2)$ is the morphism double dual to
$\varphi$, then
$$S\Psi_2S\Phi=S\Psi_1\restr{S\G(C_1)}.$$}
\end{enumerate}
Then, in some elementary extension $F^*$ of $F$, there is $C$ which is
SCF-independent from $(A,B)$ over $E$, realises $tp(C_1/A)\cup tp(C_2/B)$,
and with $S\G(C)=S_0$ (The variables for $tp(C_1/A)$ and $tp(C_2/B)$ are
identified via $\varphi$.)
\end{thm}

\prf We may assume that $F$ is sufficiently saturated; then $S\G(F)$
will also be sufficiently saturated. We work inside
$\Omega$. 
Choose $C$ realising $tp_{\SCF}(C_1/E)$, and SCF-independent from $F$
over $E$. Let $\varphi_1:C^s\to C_1^s$ and $\varphi_2:C^s\to C_2^s$ be
$E^s$-isomorphisms such that
$$\varphi_1(C)=C_1, \qquad \hbox{ and } \varphi\varphi_1=\varphi_2,
\quad(\hbox{whence}\quad 
\varphi_2(C)=C_2).$$
As $A$ is linearly disjoint from $C$ and from $C_1$ over $E$, we have that
$A^s$ is linearly disjoint from $C^s$ and from $C_1^s$ over $E^s$, and we
may therefore extend  $\varphi_1$ to an $A^s$-isomorphism
$$\varphi'_1:(AC)^s\to (AC_1)^s.$$ Similarly, we extend $\varphi_2$ to a
$B^s$-isomorphism $$\varphi'_2:(BC)^s\to (BC_2)^s.$$
Let $D_1={\varphi'_1}\inv (\acl(AC_1))$, and  $D_2={\varphi'_2}\inv
(\acl(BC_2))$. Then $D_1\subseteq (AC)^s$, $D_2\subseteq (BC)^s$.

\smallskip
Because $S\Psi_1$ and $S\Psi_2$ are elementary and $S\G(F)$ is
sufficiently saturated, there are subsystems $S_1$ and
$S_2$ of $S\G(F)$, and elementary isomorphisms
$$S\Psi'_1:S\G(\acl(AC_1))\to S_1,\quad S\Psi'_2:S\G(\acl(BC_2))\to S_2$$
extending $S\Psi_1$ and $S\Psi_2$ respectively.

\smallskip
For $i=0,1,2$ we let $L_i$ be the Galois extension of $F$ such that the restriction map $\G(F)\to
\gal(L_i/F)$ is dual to $S_i\subset S\G(F)$. 
Let $S\Phi_i'$ be the double dual of $\varphi'_i$ for $i=1,2$,  and define
$$S\Theta_i=S\Psi'_iS\Phi'_i:S\G(D_i)\to S_i,$$
and let $$\Theta_i:\gal(L_i/F)\to \G(D_i)$$
be the homeomorphism dual to $S\Theta_i$. 
We will show that there is a continuous morphism (not necessarily onto)
$$\Theta:\gal((AB)^sL_1L_2/F)\to \gal((AB)^s(AC)^s(BC)^s/ABC)$$
which induces $\Theta_i$ on $\gal(L_i/F)$ for $i=1,2$,  the identity on
$\G(\acl(AB))$, and whose image
$U$ projects onto $\G(D_1)$, $\G(D_2)$ and $\G(\acl(AB))$ (via the restriction
maps).

\smallskip
Since $A=\acl(A)$ and $B=\acl(B)$, we know that $F$ is a
regular extension of $A$ and of $B$; hence $AB$ is a regular extension
of $A$ and of $B$.
By Lemma \ref{lem6}, we may identify $\gal((AB)^s(AC)^s(BC)^s/ABC)$ with the set of triples
$(\si_1,\si_2,\si_3)\in \G(AB)\times \G(AC)\times \G(BC)$ satisfying
$$\si_1\restr{A^s}=\si_2\restr{A^s},\ \si_1\restr{B^s}=\si_3\restr{B^s},\
\si_2\restr{C^s}=\si_3\restr{C^s}.$$
For $\si\in \gal((AB)^sL_1L_2/F)$ we define
$$\Theta(\si)=(\si\restr{(AB)^s},\Theta_1(\si\restr{L_1}),
\Theta_2(\si\restr{L_2}))=:(\si_1,\si_2,\si_3).$$
We need to show that $\Theta(\si)\in\gal((AB)^s(AC)^s(BC)^s/ABC)$. Because
$\varphi'_1$ is an $A^s$-isomorphism, $S\Phi'_1$ is the identity on $S\G(A)$,
and because $S\Psi'_1$ extends $S\Psi_1$, so is $S\Psi'_1$. Hence
$S\Theta_1$ is the identity on $S\G(A)$. This shows that
$\si_1\restr{A^s}=\si_2\restr{A^s}$. Similarly,
$\si_1\restr{B^s}=\si_3\restr{B^s}$. We still need to show that
$\si_2\restr{C^s}=\si_3\restr{C^s}$. By duality, it is enough to show that
$S\Theta_1$ and $S\Theta_2$ agree on $S\G(C)$. We know that
$\varphi\varphi_1=\varphi_2$, and that $\varphi'_i$ extends
$\varphi_i$. Hence
$$S\Phi S\Phi'_1\restr{S\G(C)}=S\Phi'_2\restr{S\G(C)}.$$
We also have that $$S\Psi_1\restr{S\G(C_1)}=S\Psi_2S\Phi.$$
Hence
\begin{eqnarray*}S\Theta_1\restr{S\G(C)}&=&S\Psi_1\restr
{S\G(C_1)}S\Phi'_1\restr{S\G(C)}\cr
&=&(S\Psi_2S\Phi)(S\Phi\inv S\Phi'_2\restr{S\G(C)})\cr
&=&S\Theta_2\restr{S\G(C)},
\end{eqnarray*}
so that $\Theta$ takes its values in $\gal((AB)^s(AC)^s(BC)^s/ABC)$.
Moreover, observe that since $\varphi'_1(D_1)=\acl(AC_1)$, and by definition
of $\Theta_1$, we get that $\Theta_1$ defines a homeomorphism between
$\gal(L_1/F)$ and $\G(D_1)$. Similarly, $\Theta_2$ defines a homeomorphism
between $\gal(L_2/F)$ and $\G(D_2)$. As $\gal((AB)^sL_1L_2/F)$ projects onto
$\G(\acl(AB))$, onto $\gal(L_1/F)$ and onto $\gal(L_2/F)$ (via the restriction maps),
we get
that $U=\Theta(\gal((AB)^sL_1L_2/F))$ projects onto $\G(\acl(AB))$, onto
$\G(D_1)$ and onto $\G(D_2)$. Let $D$ be the subfield of $(AB)^s(AC)^s(BC)^s$
fixed by $U$. Then $D$ is a regular extension of $\acl(AB)$, and of $D_1$ and
$D_2$. By Theorem \ref{pac4}, there is an elementary extension $F^*$ of $F$, such that
$F^*\cap (AB)^s(AC)^s(BC)^s=D$. Note that $D_1=\acl(AC)$ and
$D_2=\acl(BC)$. 

\smallskip
To finish the proof, we need to show that $C$ realises $tp(C_1/A)\cup
tp(C_2/B)$, and that $S\G(C)=S_0$. Consider $\varphi'_1:(AC)^s\to
(AC_1)^s$. Then $\varphi'_1(C)=C_1$ and $\varphi'_1(\acl(AC))=\acl(AC_1)$
(because $F^*$ is a regular extension of $D_1$). We therefore only need to
show that the double-dual $S\Phi'_1$ is elementary in the structure
$S\G(F^*)$. By definition, $S\Phi'_1={S\Psi'_1}\inv S\Theta_1$.
By definition of $D$, the Galois extensions $L_1F^*$ and $D_1^sF^*$ are
equal, and $S\Theta_1$ is the identity on
$\gal(L_1F^*/F^*)=\gal(D_1^sF^*/F^*)$. Also, $S\Psi'_1$ is an elementary
isomorphism of $S\G(F)$, hence also of $S\G(F^*)$. This shows that
$S\Phi'_1$
is elementary, and therefore that $tp(C/A)=tp(C_1/A)$. Similarly one shows
that $tp(C/B)=tp(C_2/B)$. From the definition of $F$ one also deduces that
$F^*L_0=F^*C^s$,
which finishes the proof.

\begin{remcor}\vlabel{corfrob11} If $F$ is  Frobenius, then 
the condition on the $S\Psi_i$ can be relaxed to their being
$\call_G$-isomorphisms. The existence of $S_0$ is still required, for
trivial reasons: one could have $S\G(C_1)\subset S\G(A)$ and
$S\G(C_2)\not\subset S\G(B)$.  
\end{remcor}
\begin{defn} Let $F$ be a PAC field, let $a$, $b$, $c$ be subsets of $F$, $C=\acl(c)$, $A=\acl(c,a)$ and
$B=\acl(c,b)$. We say that $a$ and $b$ are
{\it weakly independent over $c$} if $A$ and $B$ are SCF-independent over
$C$ and $tp(S\G(A)/S\G(B))$ does not fork over $S\G(C)$.
\end{defn}

\smallskip\noindent
{\bf Remark}. Note that this notion is in general not symmetric.

\begin{thm} \vlabel{thm2} Let $F$ be a PAC field, let $E=\acl(E)\subset F$, let $a$,
$b$, $c_1$, $c_2$ be tuples of
elements of $F$, and let $A=\acl(Ea)$, $B=\acl(Eb)$ and $C_i=\acl(Ec_i)$ for
$i=1,2$. Assume that $E$ contains a $p$-basis of $F$ if the degree of
imperfection of $F$ is finite, and moreover that
\begin{enumerate}
\item[(i)]{$A\cap B=E$.}
\item[(ii)]{$A$ and $C_1$ are weakly independent over $E$, $B$ and $C_2$ are
weakly independent over $E$.}
\item[(iii)]{$tp(c_1/E)=tp(c_2/E)$.}
\item[(iv)]{by (iii), there is an $E^s$-isomorphism $\varphi:C_1^s\to
C_2^s$, sending
$c_1$ to $c_2$ and $C_1$ to $C_2$, such that the double dual map
$S\Phi:S\G(C_1)\to S\G(C_2)$ is elementary (in $S\G(F))$).
Assume that there is $S_0$ realising $tp(S\G(C_1)/S\G(A))\cup
tp(S\G(C_2)/S\G(B))$, and such that $tp(S_0/S\G(A),S\G(B))$ does not fork over
$S\G(E)$. (The variables of the two types are identified via the double dual
$S\Phi$ of $\varphi$.)}
\end{enumerate}
Then there is $c$ realising $tp(c_1/A)\cup tp(c_2/B)$, weakly independent
from $(a,b)$ over $E$.
\end{thm}

\prf Clear by Theorem \ref{thm1}.
\begin{thm}\vlabel{thm25} Let $F$ be a PAC field, and assume that $\Th(S\G(F))$
is simple. Then $\Th(F)$ satisfies the  weak independence theorem over
submodels, i.e., in the notation of \ref{thm2}, if $E\prec F$, $a$ and
$b$ are weakly independent over $E$, and
$a,b,c_1,c_2$ satisfy the hypotheses (ii) -- (iii) of \ref{thm2}, then there
is $c$ realising $tp(c_1/Ea)\cup tp(c_2/Eb)$, weakly independent from
$(a,b)$ over $E$.
\end{thm}

\prf Apply Theorem \ref{thm2}: (i) follows from the weak independence of
$a$ and $b$ over $E$, and (iv) because $S\G(F)$ satisfies the
independence theorem over models, by a result of Kim-Pillay \cite{[KP]}.

\begin{thm}\vlabel{cor25frob} Let $F$ be a Frobenius field,
$E\subset F$. Let $a,b,c_1,c_2$ be tuples in $F$, and $A=\acl(Ea)$,
$B=\acl(Eb)$, $C_1=\acl(Ec_1)$ and $C_2=\acl(EC_2)$. Assume that
\begin{enumerate}
\item[(i)]{$A\cap B=E$.} 
\item[(ii)] $a$ and $c_1$ are
SCF-independent over $E$, and $b$ and $c_2$ are SCF-independent over
$E$.
\item[(iii)]{$tp(c_1/E)=tp(c_2/E)$.}
\item[(iv)]{$\acl(S\G(A))\cap \acl(S\G(C_1))=S\G(E)$, $\acl(S\G(B))\cap \acl(S\G(C_2))=S\G(E)$.}
\end{enumerate}
Then there is $c$ realising $tp(c_1/Ea)\cup tp(c_2/Eb)$ and which is
weakly independent from $(b,c)$ over $E$.
\end{thm}

\prf By Theorem 2.4 in \cite{[C0]}, $\Th(S\G(F))$ is $\omega$-stable,
hence simple. By Proposition 4.1 in \cite{[C0]}, two subsets of $S\G(F)$ are independent over
the intersection of their algebraic closure. Apply Theorem \ref{thm25}.

\begin{remark} Condition (iv) is a little bit awkward. It is
equivalent to $S\G(A)\cap S\G(C_1)=S\G(B)\cap S\G(C_2)=S\G(E)$ in the
following two cases
\begin{itemize}
\item if $F$ is $\omega$-free, or more generally, is c-Frobenius (see
 (6.6) in \cite{[C2]} for a definition),
\item or if $E\prec F$.
\end{itemize}
\end{remark}

\prf
The first assertion follows from the fact that if $F$ is a  c-Frobenius
field, then any 
subsystem of $S\G(F)$ is algebraically closed. (In that particular case, the result
already appears in Theorem 6.4 of \cite{[C2]}.) \\
To show the second assertion, observe first that $S=S\G(E)$ is algebraically
closed: this is because $E\prec F$. Proposition 4.1 of \cite{[C0]}
tells us that if $S,S_1,S_2$ are subsystems of $S\G(F)$ with $S_1\cap
S_2=S$ and $S$ algebraically closed, then $S_1\dnfo_S S_2$, and therefore $\acl(S_1)\cap
\acl(S_2)=\acl(S)=S$.

\begin{defn} Let $n$ be an integer $>2$. A
 theory $T$
has NSOP$_n$ if for every formula $\varphi(x,y)$ (with $x$, $y$ of the
same length), if $M$ is a model of $T$, and $a_i$, $i\in \omega$, is an
infinite sequence of elements of $M$ such that $M\models
\varphi(a_i,a_j)$ whenever $i<j$, then there are $b_1,\ldots,b_n$ in
$M$ such that $M\models \varphi(b_i,b_{i+1})$ for $i=1,\ldots,n-1$,
and $M\models \varphi(b_n,b_1)$. An easy application of compactness
gives the following equivalent formulation: for all $m$ and $E\subset
M$, if $p(x,y)$ is a $2m$-type over
 such that $\bigwedge_{i\in\nat}p(x_i,x_{i+1})$ is
consistent, then $\bigwedge _{1\leq i}^{n-1}p(x_i,x_{i+1})\land p(x_n,x_1)$
is also consistent. 
\end{defn}

\begin{thm} \vlabel{thm41}
Let $F$ be a PAC field, and assume that $\Th(S\G(F))$
has NSOP$_n$ for some $n>2$. Then $\Th(F)$ satisfies NSOP$_n$.
\end{thm}

\prf If ${\rm char}(F)=p>0$ and $[F:F^p]<\infty$, we will assume that
the language contains constant symbols for elements of a $p$-basis. Assume that we have an infinite sequence $a_i$, $i\in\omega$, and a
formula $\varphi(x,y)$ such that $\varphi(a_i,a_j)$ holds whenever
$i<j$. 
We may assume that the sequence $a_i$ is indiscernible, and of length
$\aleph_1$. By stability of the theory of separably closed fields of a given
degree of imperfection, there is some $\alpha<\aleph_1$ such that for
$\beta\geq \alpha$, $tp_{\SCF}(a_\beta/a_\gamma, \gamma<\beta)$ does not fork
over $E=\acl(a_\gamma\mid \gamma<\alpha)$. Then the sequence $a_\beta$,
$\beta\geq \alpha$, is an infinite sequence of indiscernibles over
$E$.\\[0.05in] 
So, we have reduced to the case where:  we have an infinite sequence $a_i$,
$i\in\omega$, of tuples which are indiscernible and SCF-independent
over $E=\acl(E)$, and which satisfy $\varphi(a_i,a_j)$ whenever $i<j$;
moreover  $E$ contains a $p$-basis of $F$ if ${\rm char}(F)=p>0$ and
$[F:F^p]<\infty$. \\
For $i<j\in\omega$ we let $A_i=\acl(E,a_i)$ and
$K_{i,j}=\acl(A_i,A_j)$. 
Then for $i<j$, all tuples $(S\G(K_{i,j}),S\G(A_i),S\G(A_j))$ realize
the same $\call_G(S\G(E))$-type as
$(S\G(K_{0,1}),S\G(A_0),S\G(A_1))$.
We fix an $E^s$-isomorphism $\varphi_1:A_0^s\to A_1^s$, which sends
$A_0$ onto $A_1$, and denote by $S\Phi_1:S\G(A_0)\to S\G(A_1)$ the
double dual. \\[0.05in]
Since $\Th(S\G(F))$ satisfies
NSOP$_n$, there is a sufficiently saturated extension $F^*$ of $F$, and
$S_0,\ldots,S_{n-1}\subset S\G(F^*)$ such that  
 $(S_i,S_{i+1})$ and
$(S_{n-1},S_0)$ realise $tp_{\call_G}(S\G(A_0),S\G(A_1)/S\G(E))$ for $0\leq
i<n-1$. There are also $S_{i,i+1}$, $0\leq i<n-1$ and $S_{n-1,0}$ such that 
the tuples $(S_{i,j},S_i,S_j)$ realise
$tp_{\call_G}(S\G(K_{0,1}),S\G(A_0),S\G(A_1)/S\G(E))$  for $$(i,j)\in I:=\{(0,1),
 (1,2),\ldots,(n-2,n-1),(n-1,0)\}.$$  
 We fix  
$\call_G(S\G(E))$-elementary isomorphisms $S\Psi_{i,j}:
S\G(K_{0,1})\to S_{i,j}$ such that 
$$S\Psi_{i,j}(S\G(A_0))=S_i,\qquad S\Psi_{i,j}(S\G(A_1))=S_j,\qquad {S\Psi_{i,j}}\restr{S\G(A_0)}={S\Psi_{i,j}}\restr{S\G(A_1)}S\Phi_1$$ for  $i<n$, $(i,j)\in I$. \\[0.05in]
The strategy is as follows: we will use Lemma \ref{lem32} (and the
remark following it) repeatedly to
build a sequence $B_0,\ldots,B_n$ of SCF-independent realisations of
$tp(A_0/E)$, with 
$$S\G(\acl(B_iB_{i+1}))=\begin{cases}
S_{i,i+1} & \hbox{ if } i<n-1,\\
S_{n-1,0} & \hbox{ if } i=n-1,
\end{cases} 
$$
and each pair $(B_i,B_{i+1})$ realises
$tp(A_0,A_1/E)$. We 
  will then apply the amalgamation theorem \ref{thm1} to $tp(B_0/B_1)
  \cup tp(B_n/B_2\ldots
B_{n-1})$ to conclude. Let us start the construction:\\ 
By Lemma \ref{lem32} we find some $B_0$ realising $tp(A_0/E)$, an 
$E^s$-isomorphism $\psi_0:A_0^s\to B_0^s$ such that
$\psi_0(A_0)=B_0$, and the double dual  of $\psi_0$
coincides with the restriction of $S\Psi_{0,1}$ to $S\G(A_0)$. Again,
using Lemma \ref{lem32} applied to the extension $K_{0,1}$ of $A_0$, the
isomorphisms $\psi_0$ and $S\Psi_{0,1}$, we now
find some $B_1$ in $F^*$ realising $\psi_0(tp(A_1/A_0))$ and
SCF-independent from $B_0$ over $E$, an
$E^s$-isomorphism $\psi_{0,1}:(A_0A_1)^s\to (B_0B_1)^s$ extending
$\psi_0$, with $\psi_{0,1}(K_{0,1})=(B_0B_1)^s\cap F^*=:L_{0,1}$, whose
double-dual $S\G(K_{0,1})\to S\G(L_{0,1})$ coincides with $S\Psi_{0,1}$
(so that in particular $S\G(L_{0,1})=S_{0,1}$). \\[0.05in]
Induction step: At stage
$i< n$, we have found $B_0,\ldots,B_i\subset F^*$ which are  SCF-independent over $E$, and for each
$0\leq j<i$, $E^s$-isomorphisms $\psi_{j,j+1}:K_{0,1}^s\to
(B_jB_{j+1})^s$, with $\psi_{j,j+1}\varphi_1$ and $\psi_{j+1,j+2}$
agreeing on $A_0^s$ whenever $j<i-1$, such that 
$\psi_{j,j+1}(K_{0,1})=(B_jB_{j+1})^s\cap F^*$,
and the double dual of $\psi_{j,j+1}$ coincides with
$S\Psi_{j,j+1}$ on $S\G(K_{0,1})$. \\    
We now again apply Lemma \ref{lem32} to the extension $K_{0,1}/A_0$, to
the 
isomorphisms $$\psi_{i-1,i}\varphi_1: A_0^s\to B_i^s \hbox{ and } 
 S\Phi_{i,i+1}:S\G(K_{0,1})\to S_{i,i+1},$$ and find some
$B_{i+1}$ which is SCF-independent from $B_0\cdots B_i$ over $E$, an isomorphism
$\psi_{i,i+1}:(A_0 A_{1})^s\to (B_iB_{i+1})^s$ which
extends $\psi_{i-1,i}\varphi_1$,   sends $K_{0,1}$ to
$(B_iB_{i+1})^s\cap F^*$, and such that the 
 double dual of $\psi_{i,i+1}$ coincides with 
$S\Psi_{i,i+1}$. \\[0.05in]
Observe that via $S\Psi_{n-1,n}$, $S\G(B_0)=S\G(B_{n})$. By
Theorem \ref{thm1}, there is $B'_0$ realising
$$tp(B_{n}/B_2,\ldots,B_{n-1})\cup tp(B_0/B_1).$$ Then
$(B'_0,B_1,\ldots,B_{n-1})$ is our desired tuple. 

\begin{thm}\vlabel{cor42} Let $F$ be  a Frobenius field.
Then $\Th(F)$ satisfies NSOP$_3$. 
\end{thm}

\prf We know that $\Th(S\G(F))$ is $\omega$-stable by Theorem 2.4 in \cite{[C0]}, and therefore
satisfies NSOP$_3$. The result follows from Theorem \ref{thm41}.

\begin{thm}\vlabel{thm5} Let $F$ be a PAC field, let $E$, $A$, $B$
be algebraically
closed subsets of $F$, with $E$ containing a $p$-basis of $F$ if the degree
of imperfection of $F$ is finite, and assume that $A$ and $B$ are weakly
independent over $E$. Then, if $B_i$, $i\in I$, is an indiscernible
sequence of realisations of $tp(B/E)$, which is SCF-independent over $E$,
and if $p_i$ denotes the image of $tp(A/B)$ by an $E$-automorphism of $F$
sending $B$ to $B_i$, the type $\bigcup_{i\in I}p_i$ is consistent, and has a
realisation which is weakly independent from $\bigcup_{i\in I}B_i$ over
$E$. \end{thm}

\prf We may assume that $I=\nat$. Using induction and \ref{thm2}, one shows
that for every $n$, there is $A'$ realising $\bigcup_{i\leq n}p_i$, weakly
independent from $\bigcup_{i\leq n}B_i$.

\begin{remark}
 The fact that the $B_i$'s form an indiscernible sequence over
$E$ is completely unnecessary. We included this hypothesis so as to make it
look more like the usual criterion for forking. What we really prove, is
that if the $B_i$'s are SCF-independent over $E$, and for all $i\in I$,
$p_i$ is an extension of
$tp(A/E)$, having a realisation which is weakly independent from $B_i$ over
$E$, then $\bigcup_{i\in I}p_i$ has a realisation which is weakly independent
from $\bigcup_{i\in I}B_i$ over $E$.
\end{remark}

\begin{comments} Theorem \ref{cor42} also follows from results of
  Chernikov and Ramsey (Theorem 6.2 in \cite{CR}). Ramsey (\cite{R})
  shows that if $F$ is a PAC field with $S\G(F)$ NTP$_1$ or NSOP$_1$,
  then $F$ is also NTP$_1$ or NSOP$_1$. 

\end{comments} 

\begin{question}
We conclude this section with several questions. Throughout, $F$ is a
PAC field. I believe that the answer to most questions is positive. 
\begin{enumerate}
\item (Strengthening of Theorem \ref{thm41}) If $\Th(S\G(F))$ does not
  have  the strict order property, then neither does $\Th(F)$.
\item Assume that $\Th(S\G(F))$ satisfies $n$-amalgamation. Then so does
  ${\rm Th}(F)$ for the weak independence relation.
\item If $F$ is $\omega$-free, then $\Th(F)$ satisfies $n$-amalgamation
  for the strong independence relation. (Recall that $A$ and $B$ are
  strongly independent over $E$ if they are SCF-independent over $E$ and
  $\acl(AB)=\acl(A)\acl(B)$.)
\item Characterize    \thorn-forking, and whether it is
equivalent to \thorn-forking at the level of the Galois
groups. 
\item Is the theory of $\omega$-free PAC fields rosy? Superrosy? 
\item \dots 
\end{enumerate}
\end{question}



\section{More on the logic of complete systems, and codes of their
  formulas}

\para {\bf Notation}. Let $G$ be a profinite
group, $SG$ its associated system. It is convenient to consider the
equivalence relation $\sim$ associated to $\leq$: if $\alpha,\beta\in
SG$ we define $\alpha\sim \beta$ if
and only if $\alpha\leq \beta$ and $\beta\leq \alpha$. We denote by
$[\alpha]$ the $\sim$-equivalence class of $\alpha$: it comes with a
group law with graph given by $[\alpha]^3\cap P$, and for
$\beta\geq\alpha$, with a group epimorphism, denoted
$\pi_{\alpha\beta}:[\alpha]\to [\beta]$, with graph given by $C\cap
([\alpha]\times [\beta]$). The set of $\sim$-equivalence classes, equiped
with the induced partial ordering $\leq$, is a modular lattice with sup
($\lor$) and inf ($\land$). For $\alpha,\beta\in SG$, we let
$\alpha\lor\beta$ denote the identity element of $[\alpha]\lor[\beta]$,
and $\alpha\land\beta$ the identity element of
$[\alpha]\land[\beta]$. In other words: if $N_1$ is the kernel of the natural
epimorphism $G\to [\alpha]$ and $N_2$ the kernel of the natural
epimorphism $G\to [\beta]$, then $\alpha\lor\beta$ is the identity
element of $[\alpha]\lor[\beta]=G/N_1N_2$, and $\alpha\land\beta$ the identity element of
$[\alpha]\land[\beta]=G/N_1\cap N_2$. 

If $A$ is a subsystem of $SG$ and $\alpha\in SG$, we denote by
$\alpha\lor A$ the smallest (for $\leq$) element of
$\{\alpha\lor\beta\mid \beta\in A\}$. 

\para\vlabel{sg3}{\bf Action of $G$}. If $G$ is a profinite group, then $G$
acts on itself by conjugation. 
This induces an action of $G$ on $S(G)$, which respects the
$\call_G$-structure of $S(G)$. The action of an element $g$ on a given
$\sim$-equivalence class $G/N$ is then given by conjugation by $gN$, and
so does not depend on the choice of the coset representative for $gN$. 
This also defines an action of $G$ on all cartesian powers $S(G)^m$.

Let $\si$ be a tuple of
elements of $S(G)$, and $\theta(\xi,\zeta)$ an $\call_G$-formula. If
$N$ is an open subgroup of $G$ such that (the coset) $N$ is $\leq$
than all the elements of $\si$, then conjugation by the elements
of $N$ leaves the elements of the tuple $\si$ fixed, so that the
set defined by the formula $\theta(\xi,\si)$ will be invariant
under conjugation by $N$. Hence, the set defined by $\theta(\xi,\si)$
will be invariant under conjugation by $G$ if and only if, for all
$\tau\in G/N$, the sets defined by $\theta(\xi,\tau\inv\si\tau)$
and by $\theta(\xi,\si)$ coincide. Observe that in any case, the formulas
$\bigvee_{\tau\in G/N}\theta(\xi,\tau\inv\si\tau)$ and
 $\bigwedge_{\tau\in G/N}\theta(\xi,\tau\inv\si\tau)$ (with parameters
 $\si$ and $\tau$ in $G/N$) define sets which are
invariant under the action of $G$.

One can also also define an action of $G^n$ on $S(G)^{m_1}\times \cdots \times
S(G)^{m_n}$ in the natural manner.

\para\vlabel{sg6}{\bf Codes}. Let $L$ be a  Galois extension of $K$ of
degree $n$, and let 
$\si_1,\ldots,\si_m$ elements of $\gal(L/K)$. We say that a tuple of
elements of $K$ is a {\it code} for $(L,\si_1,\ldots,\si_m)$ if it
is of the form $(a,b_1,\ldots,b_m)$, where, if $a=(a_1,\ldots,a_n)$,
the polynomial $p(T)=T^n+\sum_{i=1}^{n}a_{n-i}T^i$ is the minimal monic
polynomial over $K$ of some generator $\alpha$ of $L$ over $K$, and
if $b_i=(b_{1i},\ldots,b_{ni})$, then
$\si_i(\alpha)=\sum_{j=1}^nb_{ji}\alpha^{j-1}$. 
Note that a tuple coding $(L,\si_1,\ldots,\si_m)$ will also code
$(L,\tau\inv\si_1\tau,\ldots,\tau\inv\si_m\tau)$ for any
$\tau\in\gal(L/K)$. By abuse of language, we will say that
$(L,\si,\alpha,p(T))$ is the {\it data associated to the code}
$(a,b_1,\ldots,b_m)$ ($\si=(\si_1,\ldots,\si_m)$). 
We will also say that $(a,b_1,\ldots,b_m)$ {\em codes}
$(L,\si_1,\ldots,\si_m,\alpha)$. 

Note also that the above remark shows that any orbit in $\gal(L/K)^m$ (under
the action of $\gal(L/K)$ is in fact an imaginary of $K$. 

\para\vlabel{sg7}{\bf Definable subsets of $S\G(K)$.} Let $K$
be a field. While we know (\cite{[CDM]}) that the elementary equivalence of two
fields implies the elementary equivalence of the complete systems
associated to their absolute Galois groups, the proof of this result does
not show that  subsets of  $S\G(K)$  are ``definable over
$K$''. This is easy to see: let $S\subset S\G(K)^m$ be definable over
$\gal(L/K)$, and let $\tau\in \G(K)$: then $\tau$ leaves $K$ fixed, but
its double dual sends $S$ to $\tau S\tau\inv$. One however has the following
result: 

\begin{prop} (Cherlin-van den Dries-Macintyre \cite{[CDM]}). 
Given an $\call_G$-formula
$\theta(\xi)$, which in particular says that the elements of
the tuple $ \xi$ live in the same $\sim$-equivalence class,   there is
a formula of the language of fields $\theta^*(x) $ such that for
any    field $K$, and code $a$ for an $(L,\si)$ of the appropriate sort,

  $$K\models
\theta^*(a)\iff S\G(K)\models \theta(\si).$$
\end{prop}
\noindent
The proof is easy: if $tp(b)=tp(a)$, there is an
automorphism of some ultrapower $K^\calu$ of $K$ which sends $a$ to $b$. This automorphism
extends to an automorphism of $(K^\calu)^s$, with double dual sending $\si$ to a
tuple $\si'$ coded by $ b$. Then $tp(\si)=tp(\si')$ (in $S\G(K^\calu)$
and therefore in $S\G(K)$). Hence, if
$\si$ 
satisfies $\theta(\xi)$, there is some formula $\theta_{a}( x)\in tp(a)$
which ``implies'' $\theta$, i.e., such that if $b$ satisfies $\theta_a$,
then any tuple $\si'$ coded by $b$ will satisfy $\theta$. By compactness
we get the formula $\theta^*(x)$.\\[0.15in]
The difficulties occur when one 
deals with an arbitrary
$\call_G$-formula $\theta(\xi)$, and in general one cannot hope for a
similar result. The problem is that if $\si=(\si_1,\ldots,\si_n)$ and
$a_i$ is a code for $(L_i,\si_i)$, then $a_i$ only defines $\si_i$
{\bf up to conjugation} by the elements of $\gal(L_i/K)$. Thus already a
formula of the form $\xi_i=\xi_j$ poses problem: one cannot expect to have a
formula $\theta(x_i,x_j)$ which expresses this property of {\bf all elements
coded by $x_i$    and $x_j$}. This problem can be addressed by adapting the
definition of codes, however is quite unpleasant to formulate in the general
case. Here, we will deal with a particular case.

\begin{defn}\vlabel{sg8} 
\begin{enumerate}
\item{Let $L$ be a finite Galois
extension of $K$, $\si$ a tuple of elements in $\gal(L/K)$, and
$\alpha,\beta\in L$
such that $L=K(\alpha)$. We say that $a$ {\it codes} $(L,\si,\alpha,\beta)$
if it is of the form $(b,c)$ where $b$ is a code for $(L,\si, \alpha)$ (see \ref{sg6}), and $c$ gives the coordinates of $\beta$ with respect
to the basis $\{1,\alpha,\ldots,\alpha^{[L:K]-1}\}$ of $L$ over $K$.}
\item{Let $L_1$ and $L_2$ be finite Galois extensions
of a field $K$, $\si_1$ and $\si_2$ tuples of elements in $\gal(L_1/K)$,
$\gal(L_2/K)$ respectively, and $L_0=L_1\cap L_2$. 
Let $\alpha_0,\alpha_1,\alpha_2$ be such that $L_i=K(\alpha_i)$,
$i=0,1,2$. We say that $(a_1,a_2)$ is a {\it $2$-code} for $(L_1,L_2,\si_1,\si_2,
\alpha_0,\alpha_1,\alpha_2)$ if $a_1$ is a code for
$(L_1,\si_1,\alpha_1,\alpha_0)$ and $a_2$ is a code for
$(L_2,\si_2,\alpha_2,\alpha_0)$.} 

\item{We say that $a$ is {\it a $2$-code} for $(L_1,L_2,\si_1,\si_2)$ if it
is a $2$-code for $(L_1,L_2,\si_1,\si_2,
\alpha_0,\alpha_1,\alpha_2)$ for some $\alpha_0,\alpha_1,\alpha_2$.}
\item{An $\call_G$-formula $\theta(\xi)$ is {\it codable} if it implies that the 
elements of the tuple $\xi$ are $\sim$-equivalent.}
\item{An $\call_G$-formula $\theta(\xi_1,\xi_2)$ is  {\it $2$-codable}
if it implies that the elements of the tuple $\xi_i$ are $\sim$-equivalent,
for $i=1,2$.}
\end{enumerate}
\end{defn}

\begin{remark} 
One checks easily that being a $2$-code of some $(L_1,L_2,\si_1,\si_2)$ (with
$[L_1L_2:K]\leq n$ for a fixed $n$), is an
elementary property of a tuple (again, one uses that a tuple having the same
type as a $2$-code, is a $2$-code; see also \ref{sg86}). 
One also notes that if $(a_1,a_2)$ is a $2$-code for
$(L_1,L_2,\si_1,\si_2)$, then $(a_1,a_2)$ codes exactly the tuples 
$(L_1,L_2,\rho\restr{L_1}\si_1\rho\restr{L_1}\inv,\rho\restr{L_2}\si_2\rho\restr{L_2}\inv)$
 where $\rho\in\gal(L_1L_2/K)$. Indeed, assume that $(a_1,a_2)$ codes $(L_1,L_2,\si_1,\si_2,
\alpha_0,\alpha_1,\alpha_2)$ and $(L_1,L_2,\tau_1,\tau_2,
\beta_0,\beta_1,\beta_2)$. Then $\beta_1=\rho_1(\alpha_1)$ for some
$\rho_1\in\gal(L_1/K)$, and $\beta_0=\rho_1(\alpha_0)$,
$\tau_1=\rho_1\si_1\rho_1\inv$. Similarly, $\beta_2=\rho_2(\alpha_2)$ for
some $\rho_2\in\gal(L_2/K)$, and $\beta_0=\rho_2(\alpha_0),
\tau_2=\rho_2\si_2\rho_2\inv$. This implies that $\rho_1$ and $\rho_2$ agree
on $L_0=L_1\cap L_2$, and that they can be extended to a common
$\rho\in\gal(L_1L_2/K)$. 

\end{remark}

\begin{prop}\vlabel{sg85} Let $\theta(\xi_1,\xi_2)$ be a $2$-codable
$\call_G$-formula. There is a formula
$\theta^*(x_1,x_2)$ of the language of fields, such that in any field $K$, if $(a_1,a_2)$
is a $2$-code for $(L_1,L_2,\si_1,\si_2)$, then 
$$K\models \theta^*(a_1,a_2)\iff S\G(K)\models \theta(\si_1,\si_2).$$
\end{prop}

\prf Reason as in the proof of Proposition \ref{sg7}. 

\para\vlabel{sg86}{\bf Some remarks about codes}. We saw earlier that 
being a code, or a $2$-code, is an elementary property. If one wishes 
to show this result more explicitly, one needs to work a little.

To express that $(a_1,a_2)$ is a code for $(L,\si)$ is fairly easy. 
One says first of all  that the monic polynomial $p(T)$ whose coordinates are given by 
$a_1$ is irreducible over $K$ and separable. Then, if $\alpha$ is a root of $p(T)$, 
one can interpret in $K$ the pair of fields $(K(\alpha),K)$, by 
identifying $K(\alpha)$ with $K\oplus K\alpha\oplus\cdots\oplus  K\alpha^{n-1}$ 
where $n$ is the degree of $p(T)$. One then says that $K(\alpha)$ 
contains all $n$ roots of $p(T)$, and that the tuple $a_2$ consists 
of coordinates of some of these roots (indeed, one can code an 
element $\tau$ of $\gal(K(\alpha)/K)$ by specifying the coordinates of 
$\tau(\alpha)$). 

We now want  to express the fact that $((a_1,a_2,a_3),(b_1,b_2,b_3))$ 
is a $2$-code. That $(a_1,a_2,a_3)$ is a code for some
$(L,\si,\alpha,\gamma)$ is expressible, follows from the previous 
paragraph, and similarly that $(b_1,b_2,b_3)$ is a code for some 
$(M,\tau,\beta,\delta)$. As before, one can interpret in $K$, using 
the parameters $(a_1,a_2,a_3)$ the structure $(L,K,\si,\alpha,\gamma)$, 
and similarly, using the parameters $(b_1,b_2,b_3)$, the structure 
$(M,K,\tau,\beta,\delta)$. To express that 
$((a_1,a_2,a_3),(b_1,b_2,b_3))$ is a $2$-code, we need to express the 
following:

-- that $\gamma$ and $\delta$ are conjugates over $K$.

-- that $L\cap M=K(\gamma)=K(\delta)$.

The first item is easy: in the structure 
$(K(\alpha),K,\alpha,\gamma)$, one can define the coefficients of the 
minimal (monic) polynomial $r(T)$ of $\gamma$ over $K$, and similarly one can define 
in $(M,K,\beta,\delta)$ the coefficients of the minimal polynomial 
of $\delta$ over $K$. It then suffices to say that these two minimal 
polynomials are the same, and that $K(\gamma)$ contains all roots of 
$r(T)$. 

For the second item, observe that in fact the triples 
$(L,K(\gamma),K,\alpha,\gamma)$ and $(M, K(\delta),K,\beta,\delta)$ 
are interpretable from $(a_1,a_2,a_3)$ and $(b_1,b_2,b_3)$. In 
$(L,K(\gamma),K,\alpha,\gamma)$, the coefficients of the minimal 
polynomial $q(\gamma,T)$ of $\alpha$ over $K(\gamma)$ are definable. 
It therefore suffices to say that $q(\delta,T)$ is irreducible over 
$M$, and this is expressible in the structure $(M,K,\beta,\delta)$.
[The first item gives us that $L\cap M\supseteq K(\gamma)$.  The 
irreducibility of $q(T,\delta)$ over $M$ implies that 
$[L:K(\delta)]=[LM:M]$, so that $L\cap M=K(\delta)$].

All this is done uniformly in the length of the parameters involved, 
and so gives the first-order expressibility of ``$(x,y)$ is a $2$-code''. 
Note however that the partition of the variables of the formula needs 
to be fixed, i.e., one needs to know $[L:K]=n$, $|\si|=i$, $[M:K]=m$ and 
$|\tau|=j$: if $x=(x_1,\ldots\,)$ then $a_1$ will correspond to $(x_1,\ldots,x_n)$, 
$a_2$ to $(x_{n+1},\ldots,x_{n(i+1)})$, and $a_3$ to 
$(x_{n(i+1)+1},\ldots, x_{n(i+2)})$, and similarly for the elements of 
the tuple $y$. 

\para\vlabel{sg9} {\bf An easy observation}. Let $\rho_1,\ldots,\rho_m$
enumerate an $\sim$-equivalence class of $S(G)$. Then the elements of the
subsystem of $S(G)$ generated by $\rho_1$ are in the definable closure of
$\rho_1,\ldots,\rho_m$. Indeed, each $\tau\in \langle\rho_1\rangle$ is $\geq
\rho_1$; consider the set $I(\tau)$ of indices $j$ such that
$C(\rho_j,\tau)$ holds. Because the $\rho_i$ enumerate the $\sim$-class of
$\rho_1$, the element $\tau$ is {\bf uniquely defined} by the formula
$$\bigwedge_{j\in I(\tau)}C(\rho_j,\xi)\land \bigwedge_{j\notin
I(\tau)}\neg\, C(\rho_j,\xi).$$

\begin{notation} Let $S_1,S_2$ be subsets of $S(G)$. We
denote by $tp^2(S_1/S_2)$ the set of all formulas of the form
$\theta(\xi_1,\xi_2)\in tp(S_1/S_2)$, where  the $\call_G$-formula
$\theta(\xi_1,\xi_2)$ is $2$-codable.  

Let $K$ be a field, and $A$, $B$ subfields of $K$ such that $F$ is a 
regular extension of $A$ and of $B$. We denote by $tp^*(S\G(A)/S\G(B))$ 
the set of formulas $\theta^*(X,B)$, where $\theta(\Xi_1,\Xi_2)$ is a 
$2$-codable formula, and $\theta(\Xi_1,S\G(B))\in tp^2(S\G(A)/S\G(B))$.
\end{notation} 

\begin{lem} \vlabel{sg10} Let $S_1$ and $S_2$ be subsystems of $S(G)$, and assume that 
$S_3$ satisfies $tp^2(S_2/S_1)$. Then $tp(S_3/S_1)=tp(S_2/S_1)$. 
\end{lem}

\prf By compactness, we may assume that $S_1$, $S_2$ and $S_3$ are
finite. For $i=1,2$, let $\si_i$
be an enumeration of the smallest (for $\leq$) $\sim$-equivalence class of $S_i$, and let
$\si_3$ be the subtuple of $S_3$ corresponding to $\si_2\subset S_2$. By
assumption, $tp(\si_3/\si_1)=tp(\si_2/\si_1)$, and by Observation \ref{sg9}
we get the result.  

\begin{rem}\vlabel{sg11} 
\begin{enumerate}
\item A finite disjunction of $2$-codable formulas
is not necessarily $2$-codable, but a result analogous to \ref{sg85} holds
nevertheless: \\
Let $\theta_i(\xi_i,\zeta_i)$ be $2$-codable formulas,
$i=1,\ldots,n$. Then for every field $K$, and codes $(a_i,b_i)$ for
$(L_i,M_i,\si_i,\tau_i)$, we have
$$K\models \bigvee_{i=1}^n\theta_i^*(a_i,b_i)\iff S\G(K)\models
\bigvee_{i=1}^n \theta_i(\si_i,\tau_i).$$
 \item
This result becomes false if one replaces the disjunctions by
conjunctions. However, note that a Boolean combination of $2$-codable formulas
{\bf in the same variables} is $2$-codable. 
\item From Lemma \ref{sg10}, it follows that if the field $K$ is
  $\kappa$-saturated, then so is $S\G(K)$. 
\end{enumerate}

\end{rem}

\begin{lem}\vlabel{sg12} Let $S_1$ and $S_2$ be subsystems of
some $S(G)$, and let $\Xi_1$, $\Xi_2$ enumerate the variables of $qftp(S_1)$
and $qftp(S_2)$ (the quantifier-free types). Let $\theta(\xi_1,\xi_2)$ be a
Boolean combination of $2$-codable formulas, $\xi_i\subset \Xi_i$. Then there is a
$2$-codable 
formula $\theta'(\zeta_1,\zeta_2)$, ($\zeta_i\subset \Xi_i$) such that 
$$qftp(\langle\zeta_1\rangle)\cup qftp(\langle\zeta_2\rangle)\vdash \theta(\xi_1,\xi_2)\leftrightarrow
\theta'(\zeta_1,\zeta_2).$$
\end{lem}

\prf Say that $\theta(\xi_1,\xi_2)$ is a Boolean combination of the
$2$-codable formulas $ 
\theta_i(\xi_{i1},\xi_{i2})$. Let $\zeta_1\subset \Xi_1$ enumerate a
$\sim$-equivalence class such that $qftp(S_1)\vdash \bigwedge_i (\zeta_1\leq
\xi_{i1})$, and let $\zeta_2$ be defined similarly for  $\xi_{i2}$. By
Observation 
\ref{sg9}, 
there are $2$-codable formulas $\theta'_i(\zeta_1,\zeta_2)$ such that 
$$qftp(S_1)\cup qftp(S_2)\vdash \theta_i(\xi_{i1},\xi_{i2})\leftrightarrow
\theta'_i(\zeta_1,\zeta_2).$$ 
Any  Boolean combination of the $\theta'_i(\zeta_1,\zeta_2)$ is
$2$-codable, and this gives the result (get rid of the extra variables
of $qftp(S_i)$.). 

\begin{defn}\vlabel{sg13} 
  Let $\theta(\xi)$ be a 
codable formula of $\call_G$. We define $\theta^*(x)$ to be the 
formula of the language of fields which satisfies the following 
condition, in any field $K$: \\
{\em For any tuple $a$ in $K$, $K\models \theta^*(a)$ if and only if $a$ is 
a code for some $(L,\si)$, and $S\G(K)\models \theta(\si)$.}\\[0.05in]
Similarly, if $\theta(\xi_1,\xi_2)$ is a $2$-codable formula, we let $\theta^*(x,y)$ 
be the formula of the language of fields which satisfies the 
following, for every field $K$:\\
{\em For any tuple $(a,b)$ in $K$, $K\models \theta^*(a,b)$ if and only 
if $(a,b)$ is a $2$-code for some $(L,M,\si,\tau)$ and $S\G(K)\models 
\theta(\si,\tau)$.}
\end{defn}
\noindent 
The formulas $\theta^*(x)$ and $\theta^*(x,y)$ exist, by the 
discussion above and by \ref{sg7}, \ref{sg85}. Note that $(\neg 
\theta)^*\neq \neg(\theta^*)$.

\para\vlabel{sg14}{\bf Definition of $tp^*$}. 
Let $K$ be a field, and $A$, $B$ subfields of $K$ 
such that $K\cap A^s=A$, $K\cap B^s=B$. 
We denote by $tp^*(S\G(B)/S\G(A))$ 
the set of formulas $\theta^*(X_1,A)$, where 
$\theta(\Xi_1,\Xi_2)\in tp^2(S\G(B)/S\G(A))$, and $\theta^*(X_1,X_2)$ is the 
formula of the language of fields associated to $\theta(\Xi_1,\Xi_2)$ as
in the above definition \ref{sg13}. \\
Here we need a word of explanation about the variables. The elements  
of $X_1$ correspond to an enumeration of $B$, and similarly for the 
variables of $X_2$. That $C$ satisfies $tp^*(S\G(B)/S\G(A))$ will mean 
that  we have fixed an enumeration of $C$ corresponding to the elements of 
$X_1$. From the definition of the formulas $\theta^*(X_1,X_2)$ we 
obtain the following:

\begin{rem} Assume that $C$ satisfies $tp^*(S\G(B)/S\G(A))$.
\begin{enumerate}
\item{There is an elementary $S\G(A)$-isomorphism $f:S\G(C)\to 
S\G(B)$, which respects the coding, i.e., if $a\subset A$, $b\subset 
B$ are such that $(b,a)$ is a $2$-code for some $(L_1,L_2,\si,\tau)$, and 
if $c\subset C$ is the subtuple of $C$ corresponding to $b\subset B$, 
then $(c,a)$ is a $2$-code for $(L_3,L_2,f(\si),\tau)$, where $L_3$ is defined by 
$\gal(L_1/B)=f(\gal(L_3/C))$.}

\item{If the correspondence between the elements of $B$ and the 
elements of $C$ (given by $X$) defines a field isomorphism, then $f$ is in fact 
induced by some extension of this isomorphism to an $A^s$-isomorphism 
with domain  $B^s$.}
\end{enumerate}
\end{rem}

\section{Imaginaries of PAC fields}
In this section, we will show how the type amalgamation result gives
information about imaginaries. In the later part of this chapter, we
will fix a large PAC field $F$, of characteristic $p$. If $p>0$, then
we assume that its degree of imperfection is finite, and add to the
language of rings constant symbols for elements of a $p$-basis. All
our subfields of $F$ will contain these distinguished elements.  This
has two consequences which we will use:
\begin{enumerate}
\item{The theory of separably closed fields in this expanded
language, together with axioms saying that the new constants form a
$p$-basis, is complete and eliminates imaginaries (\cite{[D]}).}
\item{If $A$ and $B$ are subfields of $F$ closed under the
$\lambda$-functions of $F$ (which give the coordinates of elements with
respect to the fixed $p$-basis), then $AB$ is also closed under the
$\lambda$-functions of $F$. This implies (4.5 in \cite{[CP]}) that
$\acl(AB)=F\cap (AB)^s$.}
\end{enumerate}
\noindent
Before starting with the description of imaginaries, we will take a
closer look at  subsets of $S\G(F)$ which are definable in $F$.

\begin{defn}\vlabel{basicim} A {\em basic
imaginary} of $F$ is a pair $( a,D)$, where $ a$ is a tuple of
elements of $F$, and $D$ is a definable subset $D$ of $S\G(F)^m$ for 
some $m$, which is stable by 
conjugation under the elements of $\G(F)$. 
\end{defn}
Here, for us, a definable or interpretable set is an imaginary element, i.e., we identify
definable/interpretable sets with their codes in the sense of $S\G(F)^{eq}$. Thus, if $L$ is a finite Galois
extension of $F$, then $\gal(L/F)$ is an imaginary, and so is $L$.

It is clear from the discussion in section 3 that basic imaginaries are 
indeed imaginaries of the field $F$. The requirement that the set $D$ 
be stable under conjugation is necessary, as elements of $\G(F)$ will 
fix elements of $F^{eq}$.

\begin{thm}\vlabel{thm31}
Let $F$ be a PAC field, of finite degree of
imperfection if the characteristic is positive, and expand the
language by adding constants for elements of a $p$-basis if necessary.
Let $e\in F^{eq}$. Then $e$ is equi-definable with a finite set of
basic imaginaries. 
\end{thm}

\prf The proof is fairly long, and proceeds with a series of steps.
We will assume that $F$ is sufficiently saturated. Let
$E=\acl^{eq}(e)\cap F$, $E_0=\dcl^{eq}(e)\cap F$. Then $E$ is a 
Galois extension of $E_0$, and every element of $\gal(E/E_0)$ lifts
to an automorphism of $F^{eq}$ fixing $e$. \\[0.05in]
If $e\in \acl^{eq}(E)$,  we are done: $e$ is coded by a tuple of
elements of $E_0$. We will therefore assume that this is not the
case, and fix a $0$-definable map $f$, and a tuple $a$ such that
$f(a)=e$. We let $A=\acl(E,a)$, and consider the set $P$ of
realisations of $tp(A/E)$. We also write $f(A)=e$.  The first step is
by now a routine argument. \\[0.05in]
We will 
consider the fundamental order on types (in the sense of $\th(F^s)$), denoted by
$\leq_{fo}$, and $\sim_{fo}$ will denote the associated equivalence
relation. We refer to chapter 13 of \cite{[P]} for the definition and properties of
the fundamental order. Recall that if $p$ and $q$ are stationary types,
then $p\sim_{fo} q$ iff they have a common non-forking extension. In our
setting, we have that if $D=\acl(D)$ and $d$ is a tuple in $F$, then
$tp_{SCF}(d/D)$ is stationary. And of course, any type in the sense of
$\th(F^s)$ over a
separably closed field  is stationary. 

\smallskip\noindent
{\bf Step 1}. There is $B\in P$, with $f(B)=e$, and which is
SCF-independent from $EA$ over $E$. \\[0.05in]
By Lemma 1.4 of \cite{[EH]}, there is $B$ realising
$tp(A/\acl^{eq}(E,e))$ such that $\acl^{eq}(E,B)\cap
\acl^{eq}(E,A)=\acl^{eq}(E,e)$, whence 
$$f(B)=e\hbox{ and }\acl(E,A)\cap
\acl(E,B)\cap F=E.\eqno{(*)}$$
Observe that because  $A$ is a regular extension of $E$, $tp_{SCF}(A/E)$ is
stationary, and so is every non-forking extension. Consider the set ${\mathcal B}$ of realisations of
$tp(A/\acl^{eq}(E,e))$ which satisfy $(*)$, and choose  $B\in {\mathcal B}$ such
that $tp_{SCF}(B/(EA)^s)$ is maximal for the fundamental order in the set
$\{tp_{SCF}(B'/(EA)^s)\mid B'\in {\mathcal B}\}$. 
Let $C$ realise
$tp(B/\acl(EA))$, SCF-independent from $B$
over $EA$. (Note that $e\in\acl(EA)$, and so $C$ realises $tp(A/\acl^{eq}(E,e))$.)
 Then $$tp_{SCF}(C/(EAB)^s)\sim_{fo}tp_{SCF}(B/(EA)^s), \hbox{ and }
f(C)=e.$$ Since $tp_{SCF}(C/(EAB)^s)\leq_{fo}tp_{SCF}(C/(EB)^s)$, we
obtain $$tp_{SCF}(B/(EA)^s) \leq_{fo}tp_{SCF}(C/(EB)^s).$$
Moreover,
$$\acl(EB)\cap \acl(EA)\cap F\subseteq \acl(EAB)\cap \acl(EC)\cap
F=\acl(EA)\cap \acl(EC)\cap F=E,$$ so that the pair $(B,C)$ satisfies $(*)$.
Since $B$ and $C$ both realise $tp(A/\acl^{eq}(E,e))$, the maximality of $tp_{SCF}(B/(EA)^s)$
among the extensions of $tp_{SCF}(A/\acl^{eq}(E,e))$ satisfying $(*)$ and the
inequality 
$$tp_{SCF}(B/(EA)^s) \leq_{fo}tp_{SCF}(C/(EB)^s)$$ imply that
$$tp_{SCF}(B/(EA)^s) \sim_{fo}tp_{SCF}(C/(EB)^s).$$ Hence
$$tp_{SCF}(C/(EAB)^s)\sim_{fo}tp_{SCF}(B/(EA)^s)\sim_{fo}tp_{SCF}(C/(EB)^s),$$
and $tp_{SCF}(C/(EAB)^s)$ does not fork over $EB$. By elimination of
imaginaries in SCF (recall that the degree of imperfection is finite,
and that $E$ contains a $p$-basis), we obtain that
$tp_{SCF}(C/EAB)$ does not fork over $\acl_{SCF}(EA)\cap
\acl_{SCF}(EB)={E}^s$, and therefore does not fork over $E$. Since
$tp(C/(EA)^s)=tp(B/(EA)^s)$, we have $tp_{SCF}(B/EA)$ does not fork over
$E$, which proves the result.

\smallskip\noindent
{\bf Step 2}. Let $B\in P$, SCF-independent from $a$ over $E$, and
with $f(B)=e$. Assume that $C\in P$ is SCF-independent from $B$ over
$E$, and that there is an $E^s$-isomorphism $\varphi:B^s\to C^s$,
whose double dual $S\Phi:S\G(B)\to S\G(C)$ is an
$\call_G(S\G(A))$-elementary map. Then $f(C)=e$. \\[0.05in]
We will apply the results of Theorem \ref{thm1}, to show that
$tp(B/A)\cup tp(C/B)$ is consistent (the identification between the
variables of $tp(B/A)$ and of $tp(C/B)$ being given by the fact that these
types extend $tp(A/E)$). 
In the notation of this theorem,
we let $\varphi=\varphi$, $S_0=S\G(C)$, $C_1=A$, $C_2=C$, $S\Psi_2$ is the identity of
$\langle S\G(C),S\G(B)\rangle$, and $S\Psi_1$ is the partial
(elementary) isomorphism on $\langle S\G(B), S\G(A)\rangle$ extending
$S\Phi$ and the identity of $S\G(A)$. \\[0.05in]
If $C'$ realises $tp(B/A)\cup tp(C/B)$, then $f(C')=f(B)=e$, and this
implies that $f(C)=e$. 

\smallskip\noindent
{\bf Step 3}. Let $B, C\in P$
with $f(B)=e$. Assume that  there is an $E^s$-isomorphism
$\varphi:B^s\to C^s$, whose double dual $S\Phi:S\G(B)\to S\G(C)$ is an
$\call_G(S\G(A))$-elementary map. Then $f(C)=e$. \\[0.05in]
By Step 1, there is  $B'\in  P$, with $f(B')=e$, and which is
SCF-independent from $A$ over $E$. As $tp(B'/A)$ has a realisation  which is
SCF-independent from $C$ over $A$,  we may assume that $B'$ is 
SCF-independent from 
$AC$ over $E$. Apply Step 2 to $A,B',C$. 

\smallskip\noindent
{\bf Step 4}. There is a $2$-codable formula $\theta(\Xi,\Upsilon)$ such
that,    
if $C\in P$, then $f(C)=e$ if and only if $S\G(C)$ satisfies
$\theta(\Xi,S\G(A))$, if and only if $C$ satisfies $\theta^*(X,A)$. \\[0.05in]
Here we need a word about the variables. If $B\in P$, then by definition we
have an 
$E$-isomorphism $\varphi:B\to A$ (which is elementary). This isomorphism
extends to an $E^s$-isomorphism $B^s\to A^s$, whose double dual is an
$S\G(E)$-isomorphism $S\G(B)\to S\G(A)$. \\[0.05in]
For each $B\in P$ such that $f(B)=e$, we know by Step 3 and by Lemma 
\ref{sg10} that $$tp(B/E)\cup tp^2(S\G(B)/S\G(A))\vdash f(X)=e.$$  
Hence there is $\theta_B(\Xi,\Upsilon)$  such that  
$\theta_B(\Xi,S\G(A))\in tp^2(S\G(B)/S\G(A))$ and 
$$tp(B/E)\cup \theta_B^*(X,A)\vdash f(X)=e.$$ By compactness, a finite
disjunction of the $\theta_B^*(X,A)$ is equivalent to $f(X)=e$ modulo
$tp(B/E)$. By Lemma \ref{sg12}, we may replace this disjunction by
$\theta^*(X,A)$ for some
$2$-codable formula $\theta(\Xi,\Upsilon)$. 

\smallskip\noindent
{\bf Step 5}. Let $B,C\in P$. Then $C$ satisfies
$\theta^*(X,B)$ if and only if $f(C)=f(B)$. \\[0.05in]
Indeed, there is an $E$-automorphism of $F$ which sends $A$ to $B$.
Then the elements of $P$ satisfying $\theta^*(X,B)$ are precisely
those satisfying $f(X)=f(B)$.

\smallskip\noindent
{\bf Step 6}. There is a  set $D$, definable over $S\G(A)$,  such
that an $E$-automorphism $\si$ of $F$ fixes $e$ if and only some (any)
extension $\tilde\si$ of $\si$ to $F^s$  leaves
$D$ invariant. \\[0.05in]
We know by Step 4 that 
$$tp(A/E)\vdash \theta^*(X,A)\iff f(X)=e.$$
Let $\si\in \aut(F/E)$ fix $e$, and $\tilde \si$ an extension of $\si$ to $F^s$. Then $\si$ induces an automorphism of
the set $P$, which leaves invariant the set of realisations of
$\theta^*(X,A)$, by step 4. Hence, (conjugation by) $\tilde \si$ leaves invariant the set $D$ of realisations
of $\theta(\Xi,S\G(A))$. I.e., $D\in \dcl^{eq}(E,e)$. \\
Conversely, let $\si\in\aut(F/E)$, and $\tilde \si$ an extension of $\si$
to $F^s$ which  leaves $D$ invariant. Then $\tilde\si$ leaves
invariant the set $P$, as well as any $\call_{S\G(E)}$-definable subset 
of $D$ which is stable by conjugation.  Hence it leaves invariant the
set of realisations of $\theta(\Xi,S\G(A))$ which are a subtuple  of some realisation of
$tp^2(S\G(A)/S\G(E))$. By Lemma \ref{lem32},  there is some $B\in P$
such that $S\G(B)$ satisfies $\theta(\Xi,S\G(A))$.  I.e., $B$ satisfies
$\theta^*(X,A)$ and $f(B)=e$. So $e\in \dcl^{eq}(E,D)$.

\smallskip\noindent
{\bf Step 7}. The result.\\[0.05in]
It follows that the imaginary $D$ and $e$ are equi-definable over $E$.
Hence there is a finite tuple $ a$ of elements of $E$ such that
they are equi-definable over $ a$. Then $e$ is equi-definable
with the set of conjugates of $( a,D)$ over $e$.

\para{\bf Remark}. This result is not totally satisfactory: it would
have been better to obtain a single basic imaginary. We will show by
an example below that it is not always possible. One can however
observe that $e$ can be squeezed between two basic imaginaries: namely
$e\in \dcl^{eq}( a,D)$, and if $b$ codes the set of
conjugates of $a$ over $e$, and $D'$ the set of conjugates of
$D$ over $e$, then $e$ is algebraic over $(b,D')$.

\bigskip\noindent
It turns out that 
the interaction between $F$ and $S\G(F)$ at the level of algebraic 
closure is very weak:  

\begin{prop}\vlabel{rem43} Let $F$ be a  PAC field. If $F$ is of  
characteristic 
$p>0$, we assume that its degree of imperfection is finite and that we 
have constant symbols for elements of a $p$-basis of $F$. 
Let 
$e=\{(a_1,D_1),\ldots,(a_n,D_n)\}\in F^{eq}$, where each $(a_i,D_i)$ is
a basic imaginary. 
\begin{enumerate}
\item {$\acl^{eq}(e)\cap F=\acl_{SCF}(a_1,\ldots,a_n)\cap F$ (=  
$\acl(a_1,\ldots,a_n)$ in the sense of the theory of $F$).}
\item{$\acl^{eq}(e)\cap 
S\G(F)=\acl^{eq}(S\G(\acl(a_1,\ldots,a_n),D_1,\ldots,D_n))\cap S\G(F)$. 
Here, as before, by $D_1,\ldots,D_n$, we mean the elements of
$S\G(F)^{eq}$ corresponding to the sets defined by $D_1,\ldots,D_n$.}
\end{enumerate}
\end{prop}

\prf (1) Using \ref{sg9},
we may assume that all the $D_i$'s are definable over $\gal(L/F)$, for
some finite Galois extension $L$ of $F$.  Let $A=\acl(a_1,\ldots,a_n)$, $b$  any finite 
tuple of elements of $F\setminus A$ such that $A(b)$ contains a code of
the extension $L$ and of the elements of $\gal(L/F)$. 
Let $B=\acl(A,b)$, and 
choose an $A^s$-automorphism $\varphi$ of $\Omega$ such that $\varphi(B)=C$ is linearly disjoint from $F$ 
over $A$. Let $\Phi:\G(B)\to \G(C)$ be the  dual of $\varphi\inv$, and 
consider the subgroup $H=\{(\si,\Phi(\si\restr{B^s}))\mid \si\in\G(F)\}$. Let 
$M$ be the subfield of $C^sF^s$ fixed by $H$. Then $M$ is a regular 
extension of $C$ and of $F$, and the restriction map $\gal(C^sF^s/M)\to 
\G(F)$ is an isomorphism. By Theorem \ref{pac4}, $F$ has an elementary extension 
$F^*$ containing $M$, regular over $M$. Then 
$tp_{F^*}(C/A)=tp_{F^*}(B/A)=tp_F(B/A)$. By definition of $H$, the 
tuple of $\gal(L/F)$ coded by $\varphi(b)$ is the same (up to conjugation) as the 
tuple coded by $b$. Moreover, as $b$ was any finite tuple of $F\setminus
A$, and $\varphi(b)\neq b$, it follows that $b\notin \acl^{eq}(e)$. \\[0.05in]
(2)  Let $S_0=\acl^{eq}(S\G(A),D_1,\ldots,D_n))\cap
S\G(F)$. Going to some sufficiently saturated extension $F^*$ of $F$ and
using again Lemma 1.4 of \cite{[EH]}, we
 find 
$S'$ realising $tp(S\G(F)/S_0)$ and such that $S'\cap S\G(F)=S_0$. (Note
that both $S'$ and $S\G(F)$ are algebraically closed). We fix some
$\call_G$-elementary map $S\Psi: S'\to S\G(F)$ which is the identity on
$S_0$. By Lemma
\ref{lem32}, there is $F_1$ realising $tp(F/A)$, and an
$A^s$-isomorphism $\psi: F_1^s\to F^s$, with double dual $S\Psi$. Since
$\psi$ is the identity on $A^s$, and $S\Psi$ is the identity on $S_0$,
it follows that $\psi(e)=e$. This shows (2).

\para\vlabel{im33}{\bf Imaginaries in complete systems of Frobenius 
fields}. Recall that a Frobenius field is a PAC field $F$, whose 
absolute Galois group $\G(F)$ has the embedding property, 
see see \ref{frob1}. The properties of $S\G(F)$ we will use are 
the following (see \cite{[C0]}, sections 2 and 4):

\begin{itemize}

\item $\th(S\G(F))$ is $\omega$-stable.

\item (Description of the types) Let $\beta$, $\gamma$ be tuples of 
elements of the equivalence class $[\beta]$, $[\gamma]$ 
respectively, let $S$ be a subsystem of $S\G(F)$, and let 
$\delta=\beta\lor S$. Then 
$tp(\beta/S)=tp(\gamma/S)$ if and only if $\gamma\lor 
S=\delta$, and there is an isomorphism 
$f:[\beta] \to [\gamma]$ such that $f(\beta)=\gamma$, and 
$\pi_{\beta,\delta}=\pi_{\gamma,\delta}f$ (i.e., $f$ induces 
the identity on $[\delta]$). 

\item If $S$ is a subsystem of $S\G(F)$, then the quantifier-free type 
of $S$ implies its type.

\item Let $A$ be a subsystem of $S\G(F)$, and $\alpha,\beta\in S\G(F)$, 
with $\alpha=\beta\lor A$. Then $\beta\in \acl(A)$ if and only if 
$\beta\in \acl(\alpha)$.
\end{itemize}

\bigskip\noindent
One cannot expect the theory of a complete system $SG$ to eliminate imaginaries,
simply because most finite groups do not eliminate imaginaries: consider
for instance $\zee/5\zee$. Then its subset $\{1,4\}$ cannot be coded by
any finite tuple of elements of $\zee/5\zee$. However, in case of
profinite groups 
with the embedding property, one obtains the next best thing: weak
elimination of imaginaries. 

\begin{thm} Let $G$ be a profinite group with the embedding property,
  $SG$ its associated system. Then ${\rm Th}(SG)$ weakly eliminate
  imaginaries. Furthermore, any imaginary is equi-definable with an
  imaginary of the form $([\alpha],\epsilon)$, where $\alpha\in SG$,
  and $\epsilon$ is an imaginary of the group $[\alpha]$. 
\end{thm}

\prf Let $D$ be a definable subset of $SG^m$, defined over some
algebraically closed 
subsystem $A$. By Observation \ref{sg9},
we may assume that if the $m$-tuple $\beta$ is in $D$, then all its
elements are $\sim$-equivalent: indeed, if  if
$(\si_1,\ldots,\si_m)\in D$, then there are only finitely many
possible isomorphism types of the system $(\langle
\si_1,\ldots,\si_m\rangle, \si_1,\ldots,\si_m)$; hence there is a
definable set $D'$ satisfying the required hypothesis,  and a
finite-to-one onto definable map $D'\to D$; then in
$SG^{eq}$, $D'\in \acl(D)$.  Since $\th(SG)$ is $\omega$-stable (by
Theorem 2.4 in \cite{[C0]}), $D$
contains only finite many types of maximal Morley rank, say
$p_i=tp(\beta_i/A)$, $i=1,\ldots,r$. Then, each $p_i$ is definable over $[\alpha_i]$,
where $\alpha_i=\beta_i\lor A$: if $B$ is  a subsystem of $SG$
containing $A$, then the unique non-forking extension of $p_i$ to $B$ is
given by $p_i\restr{[\alpha_i]}\cup \{\neg (\xi \geq \gamma)\mid
\gamma\in B, \gamma<\alpha_i\}$. I.e., $[\alpha_i]$ is the (algebraic
closure of a) canonical base for $p_i$. This shows weak elimination of
imaginaries. The last assertion follows immediately from our first reduction.

\begin{thm} 
Let $F$ be a Frobenius field, of finite invariant if the characteristic
is positive, and in that case assume that the 
language contains symbols of constants for a $p$-basis of $F$. Then
every imaginary of $F$ is equidefinable with a finite set of basic
imaginaries $(a,D)$, where furthermore $D$ is an imaginary of some group
$[\alpha]$. 
\end{thm}

\prf This is clear from the discussion above and Theorem 
\ref{thm31}. 

\para{\bf An
example}. Let $a,b,c,d$ be elements which are algebraically independent 
over $\rat$, and $\zeta$ a primitive $3$-rd root ot $1$. We fix cubic roots $\alpha_1:=\root3 \of {c+a}$,
$\beta_1:=\root3 \of {c+b}$, $\alpha_2:=\root3 \of {d+a}$, $\beta_2:=\root3 \of {d+b}$ of $(c+a)$, $(c+b)$,
$(d+a)$ and $(d+b)$ respectively. We then let
$$E_0=\rat(a,b)^{alg}, \ \ E_1=E_0(c,d, 
\alpha_2\beta_1^2, \beta_2\alpha_1^2)$$ and let $F$ be an $\omega$-free PAC field which is
a regular extension of $E_1$. Then any automorphism of $E_0$, or of
$E_1$, is elementary in the sense of ${\rm Th}(F)$. 

Let $L=F(\alpha_1,\beta_1)$, and $\si\in\gal(L/F)$ be
defined by $\si(\alpha_1)=\zeta\alpha_1$ and
$\si(\beta_1)=\zeta^2\beta_1$. Then $F(\alpha_1)=F(\beta_2)$ and
$F(\beta_1)=F(\alpha_2)$,  $\si(\alpha_2)=\zeta^2\alpha_2$, $\si(\beta_2)=\zeta\beta_2$. Consider the basic imaginaries $e_1:=(a,
(L,\si))$ and $e_2=(b,(L,\si^2))$. Note that because $\gal(L/F)$ is abelian, we do not have to
worry about conjugation. Consider the imaginary $e:=\{e_1,e_2\}$. Then
$e\in\dcl^{eq}(a,b,(L,\si))$, and letting $$f_1=(ab,a+b), \ \ 
f_2=(L,\{\si,\si^2\}),$$ we have $f_1,f_2\in\dcl^{eq}(e)$. 
 We will show that $e$ is not equidefinable with any basic imaginary
$(E,\epsilon)$. Assume by way of contradiction that
$\dcl^{eq}(e)=\dcl^{eq}(E,\epsilon)$. Then $E=\dcl^{eq}(e)\cap F$.

\smallskip\noindent
{\bf Claim}. $E=\rat(f_1)$.

We know by Proposition \ref{rem43} that $E\subset
\acl(a,b)=E_0$. Let $\rho_1$ be any automorphism of $E_1$
which fixes $c,d$ and $\zeta$, and exchanges $a$ and $b$. Then
$\rho_1(e_1)=e_2$: Indeed, $\rho_1$ extends to an automorphism
$\rho'_1$ of
$L_1:=E_1(\alpha_1,\beta_1)$, which sends $(\alpha_1,\beta_1)$ to
$(\beta_1,\alpha_1)$; one then computes
that $\rho'_1 \si {\rho'_1}\inv =\si^2$. Clearly $\rho_1(L)=L$, and so, $\rho_1(e)=e$. This being
true for any $\rho_1\in \gal(E_1/\rat(f_1,\zeta,c,d))$, we get that
$E\subseteq \rat(f_1,\zeta)$. Consider now any automorphism $\rho_2$ of
$E_1$ which fixes $a,b$, exchanges $c$ and $d$, and sends $\zeta$ to
$\zeta^2$. Then again one computes that $\rho_2(e)=e$. As $\rho_2$ moves $\zeta$
and fixes $f_1$, we obtain that $E=\rat(f_1)$. 

By Proposition \ref{rem43}(2) and because $\epsilon\in
\dcl^{eq}(a,b,(L,\si))$,  $$\acl(e)^{eq}\cap
S\G(F)=\acl^{eq}(S\G(\acl(a,b)),\epsilon)=\acl^{eq}(\gal(L/F))=\langle\gal(L/F)\rangle.$$ (The
second equality is because $S\G(\acl(a,b))=1$, and the third because any
subsystem of $S\G(F)$ is algebraically closed.)  

However, there is  an automorphism 
$\rho$ of $E_1$, which exchanges $a$ and $b$, exchanges $c$ and $d$, and is
the identity on $\zeta$. One computes that it induces the identity  on
$\gal(L/F)$, and therefore $\rho(e)\neq e$. Note that this $\rho$ is the
identity on $E$, and so $\epsilon\notin
\dcl^{eq}(E,\langle\gal(L/F)\rangle)$, which gives us the desired
contradiction.

\vskip 1cm
\noindent
Current address: \\
D\'epartement de Math\'ematiques et Applications\par\noindent
Ecole Normale Sup\'erieure\par\noindent
45 rue d'Ulm\par\noindent
75230 Paris Cedex 05\par\noindent
France\par\noindent
e-mail: {\tt zchatzid@dma.ens.fr}

\end{document}

Vieux ``Step 1''
\smallskip\noindent
{\bf Step 1}. Let $E'=\acl(E')\subset F$ contain $E$ and SCF-independent from
$A$ over $E$. There is $B\in P$, with $f(B)=e$, and which is
SCF-independent from $E'A$ over $E'$. \\[0.05in]
By Lemma 1.4 of \cite{[EH]}, there is $B$ realising
$tp(A/\acl^{eq}(E',e))$ such that $\acl^{eq}(E',B)\cap
\acl^{eq}(E',A)=\acl^{eq}(E',e)$, whence 
$$f(B)=e\hbox{ and }\acl(E',A)\cap
\acl(E',B)\cap F=E'.\eqno{(*)}$$
The first thing to note is that because $E'$ and $A$ are SCF-independent
over $E$, and $A$ is a regular extension of $E$, $tp_{SCF}(A/E')$ is
stationary. Consider the set ${\mathcal B}$ of realisations of
$tp(A/\acl^{eq}(E',e)$ which satisfy $(*)$, and choose  $B\in {\mathcal B}$ such
that $tp_{SCF}(B/(E'A)^s)$ is maximal for the fundamental order in the set
$\{tp_{SCF}(B'/(E'A)^s)\mid B'\in {\mathcal B}\}$. 
Let $C$ realise
$tp(B/\acl(E'A))$, SCF-independent from $B$
over $E'A$. 
 Then $$tp_{SCF}(C/(E'AB)^s)\sim_{fo}tp_{SCF}(B/(E'A)^s), \hbox{ and }
f(C)=e.$$ Since $tp_{SCF}(C/(E'AB)^s)\leq_{fo}tp_{SCF}(C/(E'B)^s)$, we
obtain $$tp_{SCF}(B/(E'A)^s) \leq_{fo}tp_{SCF}(C/(E'B)^s).$$
Moreover,
$$\acl(E'B)\cap \acl(E'A)\cap F\subseteq \acl(E'AB)\cap \acl(E'C)\cap
F=\acl(E'A)\cap \acl(E'C)\cap F=E,$$ so that the pair $(B,C)$ satisfies $(*)$.
Since $B$ and $C$ realise $tp(A/E')$, the maximality of $tp_{SCF}(B/(E'A)^s)$
among the extensions of $tp_{SCF}(A/E')$ satisfying $(*)$ and the
inequality \break 
$tp_{SCF}(B/(E'A)^s) \leq_{fo}tp_{SCF}(C/(E'B)^s)$ imply that
$$tp_{SCF}(B/(E'A)^s) \sim_{fo}tp_{SCF}(C/(E'B)^s).$$ Hence
$$tp_{SCF}(C/(E'AB)^s)\sim_{fo}tp_{SCF}(B/(E'A)^s)\sim_{fo}tp_{SCF}(C/(E'B)^s),$$
and $tp_{SCF}(C/(E'AB)^s)$ does not fork over $E'B$. By elimination of
imaginaries in SCF (recall that the degree of imperfection is finite,
and that $E$ contains a $p$-basis), we obtain that
$tp_{SCF}(C/E'AB)$ does not fork over $acl_{SCF}(E'A)\cap
acl_{SCF}(E'B)={E'}^s$, and therefore does not fork over $E'$. Since
$tp(C/(E'A)^s)=tp(B/(E'A)^s)$, we have $tp_{SCF}(B/E'A)$ does not fork over
$E'$, which proves the result.

\section{More results on $\omega$-free PAC fields}
In this section we will restrict our attention to a complete theory $T$
of $\omega$-free PAC fields, and will work in a large saturated model $F^*$ of $T$. If the degree of imperfection of $F^*$ is
finite, then we will add constant symbols to the language for a
$p$-basis. 

We will explore some of the notions appearing in
stability theory in the context of this theory $T$. While $T$ is not
simple, it has two pleasant notions of independence
(one the weak independence defined in \ref{thm2}, the other a  notion
stronger than forking), both satisfying the
independence theorem over algebraically closed sets. Recall also that if
the characteristic of $F$ is $0$, then  ordinary forking also
satisfies the independence theorem over algebraically closed sets.

\para{\bf Definitions and properties of forking}.  Let $C\subset A,B$ be
algebraically closed subsets of $F^*$. 
\begin{enumerate}
\item $A$ and $B$ are {\em weakly independent over} $C$ if and only if they
  are SCF-independent over $C$ and $SG(A)\cap SG(B)=SG(C)$; the last
  condition is equivalent to $F^*\cap A^sB^s=AB$, or to $A^sF^*\cap
  B^s=C^s$ (see Lemma 2.8 in \cite{[C0]}). This notion is symmetric, has
  the extension property, finite character, satisfies the independence
  theorem over algebraically closed sets. It is not transitive, nor does
  it have local character. We write $A\dnfo_C^w B$. 
\item $tp(A/B)$ does not fork over $C$ iff for any $D=\acl(D)\subset B$,
  $\acl(DA)$ and $B$ are weakly independent over $\acl(DC)$ (see
  \cite{[C0]}, theorem 3.3). This notion
  is not symmetric, nor does it have local character. It has the
  extension property, finite character. We write $A\dnfo_C^fB$. 
\item $A$ and $B$ are strongly independent over $C$ iff $A$ and $B$ are
  SCF-independent over $C$, and $\acl(AB)=AB$. This notion is symmetric,
  transitive, has finite character and the extension property, satisfies
  the independence theorem over algebraically closed sets, and in fact
  satisfies a stronger condition: if $A=\acl(A)$, and $B\supseteq A$,
  then $tp(a/A)$ has a unique strong-independent extension to
  $B$. Strong independence does not have local character. We write
  $A\dnfo_C^{s}B$. 
\item If $E,F,G\subset F^*$, we say that $E$ and $F$ are weakly independent
  [resp. strongly independent) over $G$ if the same holds of $\acl(EG)$
  and $\acl(FG)$ over $\acl(G)$. Similarly for non-forking. We then write
  $E\dnfo^*_GF$, where $*\in \{w,f, s\}$.
\end{enumerate}

\subsection{\bf Independent algebraically closed systems and $n$-amalgamation
  property}. We work in an $\omega$-free PAC-field $F$, and we will show
that their theory has the $n$-amalgamation property for the notions of
weak and strong independence.  First some definitions.

Let $n$ be an
integer, $\calp(n)$ the 
set of subsets of $\{0,1,\ldots,n-1\}$, and
$\calp^-(n)=\calp(n)\setminus \{0,1,\ldots,n-1\}$. We assume that for each
$s\in\calp^-(n)$ we are given $A_s=\acl(A_s)$, as well as
injective elementary maps $\pi_{s,t}:A_s\to A_t$ when $s\subset t$,
satisfying $\pi_{t,u}\pi_{s,t}=\pi_{s,u}$ for $s\subset t\subset u$, and
if $s,t\subset u\in\calp^-(n)$, and $v=s\cap t$, then
$\pi_{s,u}(A_s)\dnfo_{\pi_{v,u}(A_v)}\pi_{t,s}(A_t)$, and
$A_s=\acl(A_{\{i\}}, i\in s)$ (where $\dnfo$ is some independence notion). Such a family
$\{A_s,\pi_{s,t},s\subset 
t\in\calp^-(n)
\}$ is called an {\em independent system of algebraically closed sets}.

A solution for the independent system 
$\{A_s,\pi_{s,t},s\subset
t\in\calp^-(n)\}$, is some algebraically
closed $B\subset F^*$, and embeddings $\pi_s:A_s\to
A_{\{0,1,\ldots,n-1\}}$ for $s\in\calp^-(n)$ such  that
$\pi_t\pi_{s,t}=\pi_s$ whenever $s\subset t\in\calp^-(n)$, and the
structures $B_i:=\pi_{\{i\}}(A_{\{i\}})$, $i=0,\ldots,n-1$, are 
 independent over $\pi_\emptyset (A_\emptyset)$. We furthermore
require that $\pi_s$ be an elementary map for $s\in\calp^-(n)$. 

\para{\bf Theorem}. Every weakly independent system of algebraically
closed sets has a solution. Every strongly independent system of
algebraically closed sets has a solution. 

\prf  For $n=2$ there is nothing to prove. We assume the result true
for $n$ and will prove it for $n+1$. By independent, we mean weakly or
stongly independent. Let $W_0=\{s\in\calp^-(n+1)\mid
n\in s\}$. Then $W_0$ is in bijection with $\calp^-(n)$. By induction
hypothesis, there is $C=\acl(C)\subset F^*$, and elementary partial maps
$\pi_s:A_s\to B_s\subset C$ for $s\in W_0$ such that
$\pi_t\pi_{s,t}=\pi_s$ whenever $s\subset t\in W_0$; if
$C_i:=\pi_{\{i,n\}}(A_{\{i,n\}})$, then the $C_i$, $i=0,\ldots,n-1$,  are independent over
$B_n:\pi_{\{n\}}(A_{\{n\}})$; if $B_i=\pi_{\{i,n\}}(A_{\{i\}})$, then
$B_1,\ldots,B_n$ are independent over
$C_\emptyset=\pi_{\{n\}}(A_\emptyset)$. Consider now the subfields $D_0$
  and $D_1$ of $C$, defined as follows:\\[0.05in]
$D_0$ is the field composite of all $B_s:=\pi_{s\cup\{n\}}(A_s)$, where $s\in\calp^-(n)$,
and $D_1$ is the field composite of all $\pi_s(A_s)$, where $s\in
W_0$. Then $D_0\subset D_1$,  $D_1\subset (B_nD_0)^s$, $D_0$ and
$B_n$ are independent over $C_\emptyset$, and if $s\in\calp^-(n)$,
$B_{s\cup\{n\}}\subset (B_nB_s)^s$.

By induction on $n$,
there are

We will first do the strongly independent case, which is easier.

\sect{Some remarks on thorn-forking, and $1$-semisimplicity,
$2$-semisimplicity.}

In this section we will restrict our attention to a complete theory $T$
of $\omega$-free PAC fields. If the degree of imperfection of a model is
finite, then we will add constant symbols to the language for a
$p$-basis.

\para\vlabel{strdiv}{\bf Definition}. ($T=T^{eq}$) A formula $\varphi(x,b)$ {\it
strongly divides over $A$} if $b\notin \acl(A)$ and  there is an integer
$k$ such that 
$$\{\varphi(x,b')\mid tp(b/A)=tp(b'/A)\} \hbox{ is
}k\hbox{-inconsistent}.$$
In other words, if one wants to replace  imaginaries by real elements, let
$p(y)=tp(b/A)$. The
formula $\varphi(x,b)$  strongly divides over $A$ if for some $k>0$, 
$$p(y_1)\cup \cdots\cup p(y_k)\cup \{\bigwedge_{i\neq
j}\exists x\ 
\neg(\varphi(x,y_i)\leftrightarrow \varphi(x,y_j))\}\vdash \neg \exists
x\ \bigwedge_i \varphi(x,y_i).$$ 

A formula $\varphi(x,b)$ {\it \tho-divides over $A$} if it strongly divides
over $Ac$ for some tuple $c$. A formula $\varphi(x,b)$ {\it \tho-forks
over $A$} if it implies a finite disjunction of formulas, each of which
\tho-divides over $A$. A type {\tho-divides} or {\tho-forks} over $A$ if it
implies a formula which \tho-divides or \tho-forks over $A$. 

A useful criterion which we will use is the following: 

$tp(a/Ab)$ does not \tho-divide over $A$ if and only if, for every
$B\supset A$ such that $b\notin \acl(B)$, there is $a'$ realizing
$tp(a/Ab)$ and some infinite $Ba'$-indiscernible sequence containing
$b$. 

Similarly a type $q(x,a)$ over $Aa$ strongly divides over $A$ if there
is $k>0$ such that whenever $a_1,\ldots,a_k$ realize $tp(a/A)$ and are such
that the types $q(x,a_i)$ are distinct, then $\bigcup q(x,a_i)$ is
inconsistent. 

\para\vlabel{strdiv1}{\bf Lemma}. Let $F$ be a PAC field,
$E=\acl(E)\subset F$, $a$, $b$  tuples in $F$ which are SCF-independent
over $E$,  and let $A=\acl(Ea)$, $B=\acl(Eb)$. Assume that
$tp(S\G(A)/S\G(B))$ does not strongly divide over $S\G(E)$ (by which we
mean that $tp(\alpha/S\G(B))$ does not strongly divide over $S\G(E)$ for
any $\alpha\in S\G(A)$). Then $tp(a/Eb)$ does not strongly divide over
$A$.

\prf

$\varphi(x,b)$ a formula.
If $\varphi(x,b)$ strongly divides over $A$, then either $\varphi(x,b)$
SCF-forks over $A$, or \dots.

\prf 

\para\vlabel{nss}{\bf Definition}. A theory $T$ is called {\it
$1$-semisimple} if there is a notion of independence which we will
denote by $\dnfo$,  defined on subsets
of a monster model, such that:
\begin{enumerate}
\item[(i)]{(Automorphism invariance) $\dnfo$ is invariant under automorphisms.}
\item[(ii)]{(Monotonicity) $A\dnfo_CB$ if and only if $CA\dnfo_CCB$, if
and only if for all finite tuples $a\in A$ and $b\in B$, $a\dnfo_Cb$.}
\item[(iii)]{(Chain condition) For any set $A$ and  cardinal $\lambda$ there is
$\kappa$ such that whenever $(B_\alpha)_{\alpha<\kappa}$ is a family of
sets of size $<\lambda$, and $(c_\alpha)_{\alpha<\kappa}$ is a family of
finite tuples satisfying $c_\alpha\dnfo_AB_\alpha$, then there are
$\alpha\neq \beta$ and $c$ realizing $tp(c_\alpha/B_{\alpha})\cup
tp(c_\beta/B_\beta)$ with $c\dnfo_AB_\alpha B_\beta$.}
\item[(iv)]{(Existence) $A\dnfo_CC$.}
\item[(v)]{(Extension) If $c$ is finite and $B\subset B_1$, $c\dnfo_AB$, then there
is a realisation $d$ of $tp(c/AB)$ such that $d\dnfo_AB_1$.}
\end{enumerate} 
$T$ is {\it $2$-semisimple} if in addition it satisfies
\begin{enumerate}
\item[(vi)]{(Symmetry) $A\dnfo_CB$ if and only if $B\dnfo_CA$.}
\item[(vii)]{(Transitivity) If $B\subset B_1$, then $A\dnfo_CB_1$ if and
only if $A\dnfo_CB$ and $A\dnfo_BB_1$.}
\end{enumerate}

\smallskip\noindent 
{\bf Remark}. Observe that it suffices to prove (iii) for sequences
$(B_{\alpha},c_\alpha)$ which are
 indiscernible over $A$ (use Ramsey theorem). 

\para\vlabel{nss2}{\bf Theorem}. Let $F$ be a PAC field. If $\Th(S\G(F))$
 is simple, then $\Th(F)$ is $1$-semisimple. 

\prf Let $A$ and $\lambda$ be fixed, and without loss of generality,
$A\prec F$. Let $\kappa=(\lambda)^{\aleph_0})^+$, and take as
independence notion the weak independence. Then (i), (ii), (iv) and (v)
are clear (as is (vi)). We need to prove (iii).

Let
$(B_\alpha,c_\alpha)_{\alpha<\kappa}$ be indiscernible over $A$, and
such that $c_\alpha$ is weakly independent from $B_\alpha$ for every $\alpha$. 
Choose $\beta<\kappa$ such that if $E=\acl(B_\alpha\mid \alpha<\beta)$
then $B_beta$ is SCF-independent from $B_{\beta+1}$ over $E$. Such an
$\alpha \leq (\lambda^{\aleph_0})^+$ exists, by stability of SCF. Let
$d_0$ realize $tp(c_\beta/B_\beta)$ and which is weakly independent from
$\acl(EB_\beta)$ over $B_\beta$. Then $d_0$ is weakly independent from
$\acl(EB_\beta)$ over $A$. ?????

\para\vlabel{nss3} Let $F$ be an $\omega$-free PAC field, i.e., if
$F_0\equiv F$ is countable, then $S\G(F)$ is the free profinite
profinite group on $\aleph_0$ generators. They have the following two
very nice properties: 
\begin{enumerate}
\item{Let $E$ be a subfield of $F$, $a$, $b$ two tuples. Then $tp(a/E)=tp(b/E)$
if and only if there is an $E$-isomorphism $\acl(Ea)\to \acl(Eb)$ which
sends $a$ to $b$.}

\item{If $E=\acl(E)$ and $L$ is a regular extension of $E$, with
$[L:L^p]\leq [F:F^p]$ if
$char(F)=p>0$, then there is an elementary extension $F^*$ of $F$, and
an $E$-embedding $f$ of $L$ into $F^*$  such that $F^*$ is  regular
extension of $f(L)$.}
\end{enumerate}

In \cite{[C1]}, we defined the notion
of strong independence as follows: $a$ and $b$ are strongly independent
over $C$ if and only if $a$ and $b$ are SCF-independent over $C$, and
$\acl(Cab)=\acl(Ca)\acl(Cb)$. If the characteristic is $p>0$ and
$[F:F^p]<\infty$, let us add to the language constant symbols for a
$p$-basis. 

\begin{enumerate}
\item[(3)]{Then observe that for any $a$ and $C=\acl(C)\subset B$,
there is a unique extension of $tp(a/C)$ to $B$ having a realisation $a'$
which is strongly independent from $B$ over $C$: indeed, because
$C=\acl(C)$, the SCF-type $tp_{\SCF}(a/C)$ has a unique extension to $B$;
the strong independence condition requires that
$\acl(Ba')=\acl(B)\acl(Ca')$. The existence of such an extension follows
from (1) and (2) above.} 
\end{enumerate}

It then follows easily that strong independence satisfies the conditions 
(i), (ii) and (iv) -- (vii) of \ref{nss}.

\smallskip\noindent
{\bf Theorem}. The theory of an $\omega$-free PAC field is
$2$-semisimple.

\prf We take as independence notion the strong independence defined
above. We need to show \ref{nss}(iii).

Let $A$ and $\lambda$ be fixed, and without loss of generality,
$A\prec F$. Let $\kappa=2$. Let $(B_0,c_0)$ and $(B_1,c_1)$ realize the
same type over $A$, with $c_0$ strongly independent from $B_0$ over
$A$. We may assume $B_0=\acl(B_0)$. Choose $C$ realising the non-forking extension of
$tp_{\SCF}(\acl(Ac_0)/A)$ to $\acl(B_0B_1)$, and let $c\in C$ be the tuple
corresponding to $c_0$. Let $D=C\acl(B_0B_1)$. Then $D$ is a regular
extension of $\acl(B_0B_1)$, and by (2) above, we may assume that $c\in
F^*$ for some elementary extension $F^*$ of $F$, and $\acl(B_0B_1c)=D$. 
Then $c$ is strongly independent from $(B_0B_1)$ over $A$. By (3),
necessarily $c$ realises $tp(c_0/B_0)\cup tp(b_1/C_1)$.

\end{document}